%% file: 2004-26.tex
\documentclass{gtart_h}  

\input gtoutput

\lognumber{312}
\volumenumber{8}\papernumber{26}\volumeyear{2004}
\pagenumbers{969}{1012}
\received{30 March 2003}
\published{8 July 2004}
\accepted{11 June 2004}
\proposed{Shigeyuki Morita}
\seconded{Ralph Cohen, Martin Bridson}

\usepackage{amscd}
\usepackage{amssymb}
\usepackage{amsmath}
\allowdisplaybreaks
\newtheorem{thm}{Theorem}[section]
\newtheorem{lem}[thm]{Lemma}
\newtheorem{cor}[thm]{Corollary}
\newtheorem{prop}[thm]{Proposition}
\theoremstyle{definition}
\newtheorem{defn}[thm]{Definition}
\newtheorem{rem}[thm]{Remark}
\newtheorem{eg}[thm]{Example}

\def\noteq{\neq}

\def\<{\left<}
\def\>{\right>}

\def\onto{\twoheadrightarrow}

\DeclareMathOperator{\Sym}{Sym}

\DeclareMathOperator{\sgn}{sgn}

\newcommand{\field}[1]{\mathbb{#1}}
\newcommand{\ZZ}{\ensuremath{{\field{Z}}}}

\newcommand{\QQ}{\ensuremath{{\field{Q}}}}

\newcommand{\floor}[1]{\ensuremath{\left\lfloor#1\right\rfloor}}
\newcommand{\roof}[1]{\ensuremath{\left\lceil#1\right\rceil}}

\def\l{\ell}
\def\ll{\lambda}

\newcommand{\F}{\ensuremath{{\mathcal{F}}}}

\def\a{\alpha}
\def\b{\beta}
\def\g{\gamma}

\def\f{\phi}
\def\k{\kappa}

\def\s{\sigma}
\def\Sig{\Sigma}

\def\Aut{{\rm Aut}}

\def\Fat{{\F at}}
\def\wt{\widetilde}
\def\ktilde{\wt{\k}}
\def\ctilde{\wt{c}}

\def\Gam{\Gamma}

\begin{document}
\title{Increasing trees and Kontsevich cycles}
\authors{Kiyoshi Igusa\\Michael Kleber}
\address{Department of Mathematics, Brandeis
University\\Waltham, MA 02454-9110, USA}
\email{igusa@brandeis.edu, kleber@brandeis.edu}              

\begin{abstract}
It is known that the combinatorial classes in the cohomology of
the mapping class group of punctures surfaces defined by Witten
and Kontsevich are polynomials in the adjusted
Miller--Morita--Mumford classes. The leading coefficient was
computed in \cite{[I:MMM_and_Witten]}. The next coefficient was
computed in \cite{[I:GraphCoh]}. The present paper gives a
recursive formula for all of the coefficients. The main
combinatorial tool is a generating function for a new statistic on
the set of increasing trees on $2n+1$ vertices. As we already
explained in \cite{[I:GraphCoh]} this verifies all of the formulas
conjectured by Arbarello and Cornalba
\cite{[Arbarello-Cornalba:96]}. Mondello \cite{[Mondello]} has
obtained similar results using different methods.
\end{abstract}

\asciiabstract{%
It is known that the combinatorial classes in the cohomology of the
mapping class group of punctures surfaces defined by Witten and
Kontsevich are polynomials in the adjusted Miller-Morita-Mumford
classes.  The leading coefficient was computed in [Kiyoshi Igusa:
Algebr. Geom. Topol. 4 (2004) 473-520]. The next coefficient was
computed in [Kiyoshi Igusa: arXiv:math.AT/0303157, to appear in
Topology].  The present paper gives a recursive formula for all of the
coefficients. The main combinatorial tool is a generating function for
a new statistic on the set of increasing trees on 2n+1 vertices.  As
we already explained in the last paper cited this verifies all of the
formulas conjectured by Arbarello and Cornalba [J. Alg. Geom. 5 (1996)
705--749].  Mondello [arXiv:math.AT/0303207, to appear in IMRN] has
obtained similar results using different methods.}



 

\asciikeywords{Ribbon graphs, graph cohomology, mapping class group,
Sterling numbers, hypergeometric series, Miller-Morita-Mumford
classes, tautological classes}

\maketitlepage


\section*{Introduction}
\addcontentsline{toc}{section}{Introduction}

This is the last of three papers on the relationship between the
adjusted Miller--Morita--Mumford (MMM) classes $\ktilde_n$, also
known as \emph{tautological classes} (times $(-1)^{n+1}$), in the
integral cohomology of the mapping class group and certain
combinatorial classes defined by Witten and Kontsevich. In the
first paper \cite{[I:MMM_and_Witten]} we showed that these
combinatorial classes $[W_\ll^\ast]$, are polynomials in the MMM
classes and we computed the leading coefficient:
\begin{equation}\label{eq:leading coefficient for Wll}
    [W_\ll^\ast]=\prod_{i=1}^r \frac{
    ((-2)^{k_i+1}(2k_i+1)!!)^{n_i}
    }{n_i!}\ktilde_\ll+\text{ lower terms}
\end{equation}
if $\ll=k_1^{n_1}k_2^{n_2}\cdots k_r^{n_r}$ is a partition of
$\sum n_ik_i$ into $\sum n_i$ parts. Here we use the notation of
our second paper \cite{[I:GraphCoh]}
\[
    \ktilde_\ll=\prod_{i=1}^r \ktilde_{k_i}^{n_i}.
\]
The formula (\ref{eq:leading coefficient for Wll}) was conjectured
by Arbarello and Cornalba \cite{[Arbarello-Cornalba:96]} and
answers questions posed by Witten and Kontsevich
\cite{[Kontsevich:Airy]}. The introduction of
\cite{[I:MMM_and_Witten]} gives a more detailed history of the
problem.

In the next paper \cite{[I:GraphCoh]} we rephrased the theorem
(\ref{eq:leading coefficient for Wll}) above in terms of graph
cohomology using an integral version of Kontsevich's theorem that
the cohomology of the mapping class group is rationally isomorphic
to the double dual of the graph homology of connected ribbon
graphs. We also computed $a_{n,1}^{n+1}$ which is the next case of
a coefficient in the polynomial (\ref{eq:leading coefficient for
Wll}) and the dual coefficient $b_{n,1}^{n+1}$ . The notation is:
\begin{equation}\label{eq:def of coefficients allmu, bllmu}
    [W_\ll^\ast]=\sum_\mu a_\ll^\mu\ktilde_\mu,\qquad
    \ktilde_\mu^\ll=\sum_\ll b_\mu^\ll [W_\ll^\ast]
\end{equation}
where $a_\ll^\mu$ and $b_\mu^\ll$ are rational numbers.

The formula proved in \cite{[I:GraphCoh]} is
\begin{equation}\label{eq:an1 and bn1}
    a_{n,1}^{n+1}=\frac{-12a_n-(2n+5)a_{n+1}}{\Sym(n,1)},
    \qquad
    b_{n,1}^{n+1}=\frac{2n+5}{12a_n}+\frac1{a_{n+1}}
\end{equation}
where $a_n=(-2)^{n+1}(2n+1)!!$ and $\Sym(n,1)=1+\delta_{n1}$ is the
number of symmetries of $(n,1)$ (equal to 2 if $n=1$ and 1
otherwise).

The purpose of the present paper is to complete this project by
giving an algorithm for computing all of the coefficients
$a_\ll^\mu,b_\ll^\mu$ and, as an example, obtaining the following
generalization of (\ref{eq:an1 and bn1}) conjectured in
\cite{[I:GraphCoh]}.
\begin{equation}\label{eq:ank and bnk}
    a_{n,k}^{n+k}=\frac{-(2n+2k+3)a_{n+k}-a_na_k}{\Sym(n,k)},
    \qquad
    b_{n,k}^{n+k}=\frac{2n+2k+3}{a_na_k}
    +\frac1{a_{n+k}}
\end{equation}
In the meantime, Gabriele Mondello has also obtained the same
result \cite{[Mondello]}.

The contents of this paper are as follows. The first section
summarizes the definitions and results of the previous two papers.
In section~2 we study the degenerate case corresponding to degree
0 MMM class $\ktilde_0$ which is equal to the Euler characteristic
considered as a function (0--cocycle) on the space of ribbon
graphs. This is related in a simple way to the degenerate dual
Witten cycle $W_0^\ast$ which counts the number of trivalent
vertices of a ribbon graph. The formula involves Stirling numbers
of the first and second kind.

In the third section we show that the determination of the numbers
$a_\ll^\mu$ and $b_\ll^\mu$ is equivalent to the determination of
the cup product structure of the dual Kontsevich cycles. This is
more or less obvious. The coefficients in the product are not all
integers since the dual Kontsevich cycles are not integral
generators.

The coefficients $a_\ll^\mu$ are determined by the coefficients of
the inverse matrix $b_\ll^\mu$ which, by the sum of products
formula, are determined by the special cases $b_\ll^n$. Section~4
gives a formula for these coefficients $b_\ll^n$ in terms of the
category of ribbon graphs. In the next section this is reduced to
a formula involving \emph{tree polynomials}. As an example we show
in Corollary~\ref{cor:W111} that
\begin{equation}\label{eq:W111 intro}
    [W_{111}^\ast]=288\ktilde_1^3+4176\ktilde_2\ktilde_1+20736\ktilde_3
\end{equation}
This formula, together with (\ref{eq:leading coefficient for Wll})
and (\ref{eq:ank and bnk}), verifies all values of the
coefficients $a_\ll^\mu$ conjectured by Arbarello and Cornalba in
\cite{[Arbarello-Cornalba:96]}.

In Section~6 we compute the tree polynomial in the case when
almost all of the variables are equal to 1. The main application
is Section~7, where we prove the formula (\ref{eq:ank and bnk})
for $b_{r,k}^{r+k}$. The problem becomes one of finding the closed
form for a double sum of a hypergeometric term.

In Section~8 we obtain the following description of the what we
call the \emph{reduced tree polynomial}. Suppose that $T$ is an
\emph{increasing tree} with vertices $0,1,\cdots,2k$ in the sense
that, for every $0\leq j\leq2k$ the vertices $0,1,\cdots,j$ span a
connected subgraph of $T$. Then we associate to $T$ the monomial
\[
    x^T=x_0^{n_0}x_1^{n_1}\cdots x_{2k}^{n_{2k}}
\]
where $n_j$ is the number of components of $T-\{j\}$ with an even
number of vertices. The reduced tree polynomial is defined to be
\begin{equation}\label{eq:second eq for tree poly}
    \wt{T}_k(x_0,\cdots,x_{2k})=\sum_T x^T
\end{equation}
where the sum is over all increasing trees with vertices
$0,\cdots,2k$. We also show that the reduced tree polynomial
$\wt{T}_k$ is related to the tree polynomial $T_k$ of the previous
section by the formula
\[
    T_k=x_0\wt{T}_k.
\]
This tells us several things that were not obvious before. For
example, $T_k$ is a homogeneous polynomial of degree $2k+1$ with
nonnegative integer coefficients adding up to $(2k)!$. In
Section~9 we give a recursive formula for the reduced tree
polynomial. By Theorem~\ref{thm:recursive formula for blln} this
gives a recursive formula for $b_\ll^n$. By the sum of products
rule (Lemma~\ref{lem:sum of products formula for bllmu}) this
gives a formula for $b_\ll^\mu$ and thus for the $a_\ll^\mu$.
Examples are given in the last section.

The authors would like to thank Danny Ruberman for his support and
encouragement during this project. The first author is supported by
NSF Grants DMS-0204386, 0309480.

The section titles are:

\begin{enumerate}
    \item[1] Preliminaries
    \item[2] Sterling numbers and the degenerate case
    \item[3] Cup product structure of Kontsevich cycles
    \item[4] Formula for $b_\ll^n$
    \item[5] Reduction to the tree polynomial
    \item[6] First formula for $T_k$
    \item[7] A double sum
    \item[8] Reduced tree polynomial
    \item[9] Recursion for $\wt{T}_k$
    \item[10] Examples of $\wt{T}_k$
\end{enumerate}

%
%
\section{Preliminaries}

We work in the category of \emph{ribbon graphs}. These are defined
to be graphs with a designated cyclic ordering of the half edges
incident to each vertex. We consider only finite connected ribbon
graphs. We use the Conant--Vogtmann defnition \cite{[CV02]} for the
Kontsevich orientation of a connected graph. This is an ordering
up to even permutation of the set consisting of the vertices and
half-edges of the graph.

Suppose that $\Gam$ is an oriented ribbon graph and $e$ is an edge
of $\Gam$ which is not a loop (ie, the half-edges $e_1,e_2$ of
$e$ are incident to distinct vertices $v_1,v_2$. Then the graph
$\Gam/e$ obtained from $\Gam$ by collapsing $e$ to a point
$v_\ast$ has the structure of a ribbon graph and also has an
\emph{induced orientation} which is given by $v_\ast(etc.)$ if the
orientation of $\Gam$ is written as $v_1v_2e_1e_2(etc.)$. If
$\Gam$ is obtained from a trivalent graph by collapsing $n$ edges
we say that $\Gam$ has \emph{codimension} $n$.

The category of connected ribbon graphs is denoted $\Fat$. The
morphisms of this category are compositions of collapsing maps
$\Gam\to\Gam/e$ and isomorphisms. The main property of this
category is that its geometric realization is integrally homotopy
equivalent to the disjoint union of all mapping class groups
$M_g^s$ of punctured surfaces (with $s\geq1$ punctures and genus
$g$) except for the once and twice punctured sphere:
\[
    |\Fat|\simeq\coprod_{ s\geq 1,\ (s\geq 3\text{ if }g=0)
    } BM_g^s
\]
This theorem is usually attributed to Strebel \cite{[Strebel]}. A
topological proof using \emph{Outer Space} (from
\cite{[Culler-Vogtmann-86]}) can be found in \cite{[I:BookOne]}.

By a theorem of Kontsevich proved in \cite{[CV02]} and refined in
\cite{[I:GraphCoh]}, the cohomology of $\Fat$ (or equivalently,
$M_g^s$) is rationally isomorphic to the cohomology of the
associative graph cohomology complex. We work in the \emph{integer
subcomplex} of the rational associative graph cohomology complex
generated by the cochains
\[
    \<\Gam\>:=|\Aut(\Gam)|[\Gam]^\ast
\]
This is a $\ZZ$--augmented complex of free abelian groups which can
be described as follows.

\begin{defn}
For all $n\geq0$ let $G_n^\ZZ$ be the free abelian group generated
by all isomorphism classes $\<\Gam\>$ of oriented connected ribbon
graphs $\Gam$ of codimension $n$ without orientation reversing
automorphisms modulo the relation $\<-\Gam\>=-\<\Gam\>$. For
$n\geq1$ let $d\co G_n^\ZZ\to G_{n-1}^\ZZ$ be given by
\[
    d\<\Gam\>=\sum \<\Gam_i\>
\]
where the sum is over all isomorphism classes of oriented ribbon
graphs $\Gam_i$ over $\Gam$ with one extra edge $e_i$ so that
$\Gam\cong\Gam/e_i$ with the induced orientation.
\end{defn}

\begin{thm}[Kontsevich \cite{[CV02]}]\label{thm:Kontsevic}
$H^\ast(\coprod BM_g^s;\QQ)\cong H^\ast(G_\ast^\ZZ;\QQ)$.
\end{thm}

The refinement of this theorem proved in \cite{[I:GraphCoh]} is:

\begin{thm}
This rational equivalence is induced by an augmented integral
chain map
\[
    \f\co C_\ast(\Fat)\to G_\ast^\ZZ
\]
where $C_\ast(\Fat)$ is the cellular chain complex of the nerve of
$\Fat$.
\end{thm}

If $\ll=1^{r_1}2^{r_2}\cdots$ is a partition of $n=\sum ir_i$, the
\emph{dual Kontsevich cycles} $W_\ll^\ast$ is the integral $2n$
cocycle on the integral cohomology complex $G_\ast^\ZZ$ given as
follows:
\[
    W_\ll^\ast(\<\Gam\>)=o(\Gam)=\pm1
\]
if $\Gam$ is an oriented ribbon graph of codimension $2n$ having
exactly $r_i$ vertices of valence $2i+3$ and no even valence
vertices. The sign is $+$ if $\Gam$ has the \emph{natural
orientation} (given by taking each vertex followed by the incident
half edges in cyclic order) and $-$ is not. This set of ribbon
graphs is denoted $W_\ll$ and called the \emph{Kontsevich cycle}.
If $\Gam$ is not in $W_\ll$ then $W_\ll^\ast(\<\Gam\>)=0$.

Recall that the \emph{Miller--Morita--Mumford class} $\k_n\in
H^{2n}(BM_g,\ZZ)$ is defined topologically (\cite{[Miller86:MMM]},
\cite{[Morita84]}) as the image under the transfer
\[
    p_\ast\co H^{2n+2}(E)\to H^{2n}(BM_g)
\]
of the $n-1^{\rm st}$ power $e^n$ of the Euler class $e\in H^2(E)$ of
the vertical tangent bundle of the universal surface bundle over
$BM_g$ with fiber an oriented surface $\Sig_g$ of genus $g$. If we
pull this surface bundle back to the space $B=BM_g^s$ which maps
to $BM_g$, we get $s$ points in each fiber forming an $s$--fold
covering space $\wt{B}$ over $B$. The \emph{adjusted} or
\emph{punctured} Miller--Morita--Mumford class is given by
\[
    \ktilde_n=\k_n-p_\ast(c^{n})
\]
where $c\in H^2(\wt{B})$ is the Euler class of the vertical
tangent bundle of $E$ pulled back to $\wt{B}$. (See
\cite{[I:ComplexTorsion]} for more details about this construction
and its relationship to higher Franz--Reidemeister torsion.)
Arbarello and Cornalba \cite{[Arbarello-Cornalba:96]} showed that
these are the correct versions of the MMM classes which should be
compared to the combinatorial classes of Witten and Kontsevich.

In \cite{[I:MMM_and_Witten]} it was shown that the adjusted MMM
classes are represented by the \emph{cyclic set cocycle}
$c^n_\Fat$ adjusted by a factor of $-2$:
\[
    \ktilde_n=-\frac12[c_\Fat^n].
\]
Therefore, $\ktilde_n$ is represented by the \emph{adjusted cyclic
set cocycle}
\[
    \ctilde_n=-\frac12c_\Fat^n.
\]
This cocycle can be defined as follows. Take any $2n$--simplex
\[
    \Gam_\ast\co \Gam_0\to\Gam_1\to\cdots\to\Gam_{2n}
\]
in the category of ribbon graphs. Then
\[
    \ctilde_n(\Gam_\ast)=-\frac12\sum_v
    m(v)\sum\frac{\sgn(a_0,a_1,\cdots,a_{2n})}
    {|C_0|\cdot|C_1|\cdots|C_{2n}|}
\]
where the first sum is over all vertices $v$ of $\Gam_0$, $m(v)$
is the valence of $v$ minus 2, and the second sum is over all
choices of angles $a_i$ of the vertex $v_i$ which is the image of
$v$ in $\Gam_i$. The denominator has the sizes $|C_i|$ are the
sets $C_i$ of angles about $v_i$ (so $a_i\in C_i$ for each $i$).
The sign is the sign of the permutation of the images of $a_i$ in
the final set $C_{2n}$. When these angles are not distinct, the
sign is zero and, more generally, the sign sum is equal to the
partial sum given by choosing each $a_i$ in the complement of the
image of $C_{i-1}$ in $C_i$. For more details, see
\cite{[I:MMM_and_Witten]}.

The relationship between the adjusted MMM classes $\ktilde_n$ and
the \emph{dual Witten cycles} $[W_n^\ast]$ is given
(\cite{[I:MMM_and_Witten]}) by
\[
    [W_n^\ast]=a_n\ktilde_n,\qquad \ktilde_n=b_n[W_n^\ast]
\]
where
\[
    a_n=\frac1{b_n}=(-2)^{n+1}(2n+1)!!
\]
To compute the other coefficients in (\ref{eq:def of coefficients
allmu, bllmu}) we need the following formula proved in
\cite{[I:GraphCoh]}, Lemma~3.15.
\begin{lem}[Sum of
products rule]\label{lem:sum of products formula for bllmu} If
$\ll=(\l_i,\cdots,\l_r)$ is a partition of $n$ into $r$ parts and
$\mu=(m_1,\cdots,m_s)$ is a partition of the same number $n$ into
$s$ parts then the coefficient $b_\ll^\mu$ in equation
(\ref{eq:def of coefficients allmu, bllmu}) is equal to the sum
\[
    b_\ll^\mu=\sum_f\prod_{j=1}^s b_{\ll_{\pi(j)}}^{m_j}
\]
over all epimorphisms
\[
    f\co \{1,\cdots,r\}\onto\{1,\cdots,s\}
\]
having the property that the sum of the numbers $\l_i$ over all
$i\in \pi(j)=f^{-1}(j)$ is equal to $m_j$ of the product over all
$1\leq j\leq s$ of the coefficient $b_{\ll_{\pi(j)}}^{m_j}$ where
${\ll_{\pi(j)}}$ is the partition of $m_j$ given by the numbers
$\l_i$ for $i\in \pi(j)$.
\end{lem}
By this formula it suffices to compute the numbers $b_\ll^m$.

%
%

\section{Sterling numbers and the degenerate case}

We start with an examination of the degenerate case
$W_{0^n}^\ast$. These are polynomials in the $0^{\rm th}$ adjusted
cyclic set cocycle $\ctilde_0$, equal to the $0^{\rm th}$ (topological)
Miller--Morita--Mumford class $\ktilde_0$, which is the Euler
characteristic. If $\Gam$ is trivalent with the natural
orientation, then
\[
    \ctilde\<\Gam\>=\chi(\Gam)=\frac{v}{-2}
\]
where $v$ is the number of vertices of $\Gam$. (In general we need
to count the number of vertices with \emph{multiplicity}, ie,
valence minus $2$.)

We interpret the $0$'s in $W_{0^n}^\ast$ as counting the number of
vertices with multiplicity:
\begin{equation}\label{eq:W0n = S1 c0i}
    W_{0^n}^\ast\<\Gam\>=\binom{v}{n}=\binom{-2\ctilde_0}{n}=\frac1{n!}\sum_{i=0}^nS_1(n,i)(-2\ctilde_0)^i
\end{equation}
where $S_1(n,i)$ is the Stirling number of the first kind. This
can be solved for the $\ctilde_0^i$ to give:
\begin{equation}\label{eq:c0m = S2 W0n}
    \ctilde_0^m=\frac1{(-2)^m}\sum_{n=0}^m n!S_2(m,n)W_{0^n}^\ast
\end{equation}
where $S_2(m,n)$ are the Stirling numbers of the second kind.

In the notation of \cite{[I:GraphCoh]}, this is
\[
    \ctilde_0^m=\sum_{n=0}^m b_{0^m}^{0^n}W_{0^n}^\ast
\]
where
\begin{equation}\label{eq:b0m0n = S2}
    b_{0^m}^{0^n}=\frac{n!S_2(m,n)}{(-2)^m}.
\end{equation}
This is consistent with the formula
\[
    b_{0^m}^{0^n}=\sum_f \prod_{j=1}^n b_{0^{m_j}}^0=\sum_f
    \prod_{j=1}^n\frac1{(-2)^{m_j}}
\]
where the sum is taken over all surjective mappings
\[
    f\co \{1,2,\cdots,m\}\twoheadrightarrow
    \{1,2,\cdots,n\}
\]
with $m_j$ being the number of elements in $\pi(j)=f^{-1}(j)$.
Since there are $n!S_2(m,n)$ such mappings $f$, this agrees with
(\ref{eq:b0m0n = S2}).

Assume for a moment that the sum of products formula
(Lemma~\ref{lem:sum of products formula for bllmu}) holds more
generally for all partitions with $0$'s. Thus, if
$\mu=(\mu_1,\mu_2,\cdots,\mu_r)$ and
$\ll=(\ll_1,\ll_2,\cdots,\ll_s)$ are partitions of the same number
$n$ then we have the following which we take as a definition. (It
agrees with the previously defined terms $b_\ll^\mu$ when
$p=q=0$.)
\begin{equation}\label{eq:bmu0m ll0n}
    b_{\ll0^p}^{\mu0^q}:=\sum_f b_{\ll_{\pi(1)}}^{\mu_1}\cdots
    b_{\ll_{\pi(r)}}^{\mu_r}
    b_{\ll_{\pi(r+1)}}^0\cdots b_{\ll_{\pi(r+q)}}^0
\end{equation}
where the sum is over all surjective mappings
\[
    f\co \{1,2,\cdots,s+p\}\twoheadrightarrow
    \{1,2,\cdots,r+q\}
\]
having the property that the sum of the parts $\ll_j$ of $\ll$ for
$j\in\pi(i)=f^{-1}(i)$ is equal to $\mu_i$:
\[
    \mu_i=\sum_{j\in\pi(i)}\ll_j
\]
where $\ll_j=0$ for $i>s$ and $\mu_i=0$ for $i>r$. When the
superscript of $b$ is $0$ the subscript must be $0^m$ for some
$m\geq1$ and we have
\[
    b_{0^m}^0=\frac1{(-2)^m}.
\]
If the superscript is $\mu_i\noteq0$ then the subscript is a
partition of $\mu_i$, say $\nu$, plus any number of $0$'s. We
define
\[
    b_{\nu0^m}^{\mu_i}:=\frac{(2\mu_i+1)^m}{(-2)^m}b_{\nu}^{\mu_i}.
\]
This makes sense since it is supposed to be the contribution of a
vertex of valence $2\mu_i+3$ to the cup product
\[
    \ktilde_{\nu0^m}=\ktilde_\nu\ktilde_0^m
\]
But each $\ktilde_0$ is given by
\[
    \frac{v}{-2}=\frac{2\mu_i+1}{-2}.
\]
Putting these together in (\ref{eq:bmu0m ll0n}) we get the
following.

\begin{prop}\label{prop:bll0p mu0q =S2}
\[
    b_{\ll0^p}^{\mu0^q}=\sum_{m=0}^{p-q}\binom{p}{m}q!S_2(p-m,q)\frac{(2n+r)^m}{(-2)^p}b_\ll^\mu
\]
\end{prop}

We claim that these are the coefficients which convert monomials
in the adjusted Miller--Morita--Mumford classes into linear
combinations of dual Kontsevich cycles with $0$'s.

\begin{defn}\label{def:Kontsevich cycles with 0s} Let
$\ll=1^{n_1}2^{n_2}\cdots$ be a partition of $n=\sum in_i$ into
$r=\sum n_i$ parts. We define the \emph{degenerate Kontsevich
cycles} $W_{\ll0^m}^\ast$ to be the integer cocycle of degree $2n$
on the integer subcomplex of associative graph cohomology given by
\[
    W_{\ll0^m}^\ast\<\Gam\>=o(\Gam)\binom{n_0}{m}
\]
provided that $\Gam$ is a connected oriented ribbon graph having
exactly $n_i$ vertices of valence $2i+1$ for all $i\geq0$ and no
vertices of even valence. The orientation $o(\Gam)$ is $\pm1$
depending on whether or not the orientation of $\Gam$ is the
natural one.
\end{defn}

It is easy to express $W_{\ll0^m}^\ast$ in terms of the Euler
characteristic
\[
    \chi=\ctilde_0=\frac{n_0+2n+r}{-2}
\]
and the nondegenerate Kontsevich cycle $W_\ll^\ast$:
\[
    W_{\ll0^m}^\ast=\frac1{m!}\sum_{j=0}^m
    S_1(m,j)(-2\ctilde_0-2n-r)^jW_\ll^\ast
\]
\[
    =\frac1{m!}\sum_{0\leq i\leq j\leq m}
    S_1(m,j)\binom{j}{i}(-2n-r)^{j-i}(-2\ctilde_0)^iW_\ll^\ast
\]
Passing to cohomology classes, this can be written as follows.

\begin{thm} The degenerate Kontsevich cycles are related to the
adjusted Miller--Morita--Mumford classes by
\[
    [W_{\ll0^m}^\ast]=\sum_{\mu,i}
    a_{\ll0^m}^{\mu0^i}\ktilde_\mu\ktilde_0^i
\]
and
\[
    \ktilde_\ll\ktilde_0^p=
    \sum_{\mu,q}b_{\ll0^p}^{\mu0^q}[W_{\mu0^q}^\ast]
\]
where
\[
    a_{\ll0^m}^{\mu0^i}=\frac1{m!}\sum_{j=i}^m
    S_1(m,j)\binom{j}{i}(-2n-r)^{j-i}(-2)^i a_\ll^\mu
\]
and $b_{\ll0^p}^{\mu0^q}$, defined by (\ref{eq:bmu0m ll0n}), is
given by Proposition \ref{prop:bll0p mu0q =S2}.
\end{thm}

\begin{proof}
Using the duality between the first and second Stirling numbers it
is easy to see that the matrices with coefficients
$b_{\ll0^p}^{\mu0^q}$, $a_{\mu0^q}^{\ll0^p}$ are inverse to each
other.
\end{proof}

%
%

\section{Cup product structure of Kontsevich cycles}

Using Kontsevich's theorem (\ref{thm:Kontsevic}) the rational
cohomology of $G_\ast^\ZZ$ inherits a ring structure.

\begin{thm}The determination of the conversion coefficients $a_\ll^\mu$
and $b_\mu^\ll$ is equivalent to finding the coefficients
$m_{\ll\mu}^\nu$ giving the cup product of the Kontsevich
cocycles:
\[
    [W_\ll^\ast]\cup[W_\mu^\ast]=\sum_{\nu}m_{\ll\mu}^\nu
    [W_\nu^\ast]\in H^\ast(G_\ast;\QQ)
\]
\end{thm}

\begin{rem}Note that rational numbers $m_{\ll\mu}^\nu$ are well-defined
since $[W_\ll^\ast]$ are linearly independent over $\QQ$ and span
the same vector subspace as the monomials in the adjusted
Miller--Morita--Mumford classes $\ktilde_\ll$. We also note that
these numbers are not all integers. The simplest example is
\[
    [W_1^\ast]\cup[W_1^\ast]=2[W_{1,1}]+\frac{29}5[W_2^\ast]
\]
which follows from the equations:
\[
    [W_1^\ast]=a_1\ktilde_1=12\ktilde_1
\]
\[
    \ktilde_1^2=2(b_1)^2[W_{1,1}^\ast]+b_{1,1}^2[W_2^\ast]
\]
\[
    =\frac2{144}[W_{1,1}^\ast]+\left(
\frac7{144}-\frac1{120}
    \right)[W_2^\ast]
\]
\end{rem}

\begin{proof}In one direction this is clear. If we know the numbers
$a_\ll^\mu$ and $b_\mu^\ll$ then we can convert $[W_\ll^\ast]=\sum
a_\ll^\a\ktilde_\a$ and $[W_\mu^\ast]=\sum a_\mu^\b\ktilde_\b$,
multiply and convert back. Thus,
\begin{equation}\label{eq:multiplication coefficients}
    m_{\ll\mu}^\nu=\sum_{\a,\b}a_\ll^\a a_\mu^\b b_{\a\b}^\nu.
\end{equation}
The other direction is also easy. Suppose we know the numbers
$m_{\ll\mu}^\nu$ and we want to find $a_\ll^\mu$, $b_\ll^\mu$. We
proceed by induction on the number of parts of $\ll$. When $\ll=n$
is a partition of $n$ with one part, then $\mu$ must also be equal
to $n$ since $\mu$ cannot have more parts than $\ll$. But we know
these numbers:
\[
    a_n^n=\frac1{b_n^n}=(-2)^{n+1}(2n+1)!!
\]
Suppose by induction that we know $a_\ll^\mu$, $b_\ll^\mu$ for all
partitions $\ll$ with $r$ or fewer parts. Then setting $\mu=n$ in
(\ref{eq:multiplication coefficients}) there will be only one term
on the right hand side (when $\a=\ll$ and $\b=\mu=n$) which is
unknown. This gives $b_{\ll n}^\nu$. Taking the inverse matrix we
also get all $a_{\ll n}^\nu$.
\end{proof}

%
%

\section{Formula for $b_\ll^n$}

Using the sum of products rule (Lemma~\ref{lem:sum of products
formula for bllmu}), the calculation of the numbers $b_\ll^\mu$ is
reduced to the case when $\mu=n$ is a partition of $n$ with one
part. If $\ll=(\ll_1,\cdots,\ll_r)$ is a partition of $n$ into $r$
parts then the number $b_\ll^n$ is given by
\[
    b_\ll^n =(-1)^n\ctilde_\ll D(\Gam)=(-1)^n(\ctilde_{\ll_1}\cup
    \cdots\cup\ctilde_{\ll_r}) D(\Gam)
\]
where $\Gam$ is a ribbon graph with natural orientation having one
vertex of valence $2n+3$ and all other vertices trivalent and
$D(\Gam)$ is any dual cell of $\Gam$.
\[
    D(\Gam)=\sum_{\Gam_\ast} o(\Gam_\ast)
    (\Gam_0\to\cdots\to\Gam_{2n}=\Gam)
\]
where the sum is over all sequences of morphisms over $\Gam$
between representatives $\Gam_i$ of the isomorphism classes of
ribbon graphs over $\Gam$ and $o(\Gam_\ast)=\pm1$ is positive iff
the natural orientations of $\Gam=\Gam_{2n}$ agrees with the
orientation induced from the natural orientation of the trivalent
graph $\Gam_0$ by the collapsing morphisms in the sequence
$\Gam_\ast=(\Gam_0\to\cdots\to\Gam_{2n}=\Gam)$ which we abbreviate
as $(\Gam_0,\cdots,\Gam_{2n})$.

Combining these we get
\[
    b_\ll^n=(-1)^n\sum_{\Gam_\ast} o(\Gam_\ast)
    \ctilde_{\ll_1}(\Gam_0,\cdots,\Gam_{2\ll_1})
    \cdots
    \ctilde_{\ll_r}(\Gam_{2n-2\ll_r},\cdots,\Gam_{2n})
\]
We use the notation
\[
    \ll[i]=\ll_1+\ll_2+\cdots+\ll_i
\]
(with $\ll[0]=0$ and $\ll[r]=n$). Then the $i^{\rm th}$ factor in the
expression for $b_\ll^n$ is
\begin{equation}\label{eq:ith factor of bll n without sign}
    \ctilde_{\ll_i}(\Gam_{2\ll[i-1]},\Gam_{2\ll[i-1]+1},
    \cdots,\Gam_{2\ll[i]})
\end{equation}
We will factor the sign terms $(-1)^n$ and $o(\Gam_\ast)$ into $r$
factors and associate each factor to one of the factors
(\ref{eq:ith factor of bll n without sign}).

First, we note that the graphs $\Gam_{\ll[i]}$ must all be odd
valent in the sense that they have no even valent vertices. If not
then one of the $\ctilde_{\ll_i}$ factors (\ref{eq:ith factor of
bll n without sign}) would be zero. Consequently, the orientation
term $o(\Gam_\ast)$ can be factored as:
\[
    o(\Gam_\ast)=\prod_{i-1}^r
    o(\Gam_{2\ll[i-1]},\cdots,\Gam_{2\ll[i]}).
\]
The sign $(-1)^n$ also factors:
\[
    (-1)^n=\prod (-1)^{\ll_i}.
\]
So, $b_\ll^n$ is a sum of products. Each product has $r$ factors
where the $i^{\rm th}$ factor has the form
\begin{equation}\label{eq:simple factors for bll n}
    (-1)^{\ll_i}o(\Gam_{2\ll[i-1]},\cdots,\Gam_{2\ll[i]})
    \ctilde_{\ll_i}(\Gam_{2\ll[i-1]},\cdots,\Gam_{2\ll[i]})
\end{equation}
which we abbreviate as $(-1)^{\ll_i}o(\Gam_\ast^i)
\ctilde_{\ll_i}(\Gam_\ast^i)$. But, the graphs
$\Gam_{2\ll[i-1]+j}$ for $1\leq j<\ll_i$ occur only in the $i^{\rm th}$
factor (\ref{eq:simple factors for bll n}). Thus, we have the
following.

\begin{lem}\label{lem:bll n as a sum of products of sums}
$b_\ll^n$ can be expressed as a sum of products of sums:
\begin{equation}\label{eq:bll n as a sum of products of sums}
    b_\ll^n=\sum_{(\Gam_0,\Gam_{2\ll[1]},\cdots,\Gam_{2\ll[r]})}
    \prod_{i=1}^n
    \sum_{\Gam_\ast^i}
    (-1)^{\ll_i}o(\Gam_\ast^i)
    \ctilde_{\ll_i}(\Gam_\ast^i)
\end{equation}
The first summation is over all sequences $\Gam_{2\ll[i]}$,
$i=0,\cdots,r$ of (representatives of) isomorphism classes of odd
valent graphs over $\Gam=\Gam_{2n}=\Gam_{2\ll[r]}$ and the second
sum is over all sequences of morphisms
\[
    \Gam_\ast^i=(\Gam_{2\ll[i-1]}\to\Gam_{2\ll[i-1]+1}\to\cdots
    \to\Gam_{2\ll[i]})
\]
where each $\Gam_j$ is from a fixed set of representatives from
the set of isomorphism classes of oriented ribbon graphs over
$\Gam$.
\end{lem}

Now we examine the possibilities for the graphs $\Gam_{2\ll[i]}$.
Since $\Gam_0$ has codimension $0$ it must be trivalent. In order
for the first factor in (\ref{eq:bll n as a sum of products of
sums}) to be nonzero, we must have that $\Gam_{2\ll[1]}$ is
trivalent except at one vertex of valence $2\ll_1+3$. More
generally, we have the following.

\begin{lem}\label{lem:valences of vertices of Gam ll(i)}
Suppose that the nontrivalent vertices of $\Gam_{2\ll[i]}$ have
valences $2n_{i1}+3,\cdots,2n_{ik_i}+3$. Then, in order for the
corresponding terms in (\ref{eq:bll n as a sum of products of
sums}) to be nonzero we must have the following.
\begin{enumerate}
    \item $\ll(i)=n_{i1}+\cdots+n_{ik_i}$
    \item For each $i\geq1$ and each $j<k_i$ there is an index
    $\f(j)$ so that $n_{i-1,\f(j)}=-n_{ij}$ and $\f$ is an injective function.
\end{enumerate}
\end{lem}

\begin{proof}
In order for the term
$\ctilde_{\ll_i}(\Gam_{2\ll[i-1]}\to\cdots\to\Gam_{2\ll[i]})$ to
be nonzero, the inverse images in $\Gam_{2\ll[i-1]}$ of the
vertices of $\Gam_{2\ll[i]}$ must all be vertices (necessarily
with the same valence) with only one exception. The exceptional
vertex must have valence at least $2\ll_i+3$ and its inverse image
must be a tree in $\Gam_{2\ll[i-1]}$ with that many leaves.
\end{proof}

Next, we look at the factors in (\ref{eq:bll n as a sum of
products of sums}) for $i=1,\cdots,r$. The first factor is easy to
compute:
\[
    \sum_{\Gam^1_\ast}
    (-1)^{\ll_1}o(\Gam^1_\ast)\ctilde_{\ll_i}(\Gam^1_\ast)
    =
    b_{\ll_1}=\frac1{(-2)^{\ll_1+1}(2\ll_1+1)!!}
\]
The last factor ($i=r$) is more difficult. It is also
\emph{universal} in the sense that, if we can compute the last
factor, we can compute all the factors. We make this statement
more precise using the tree polynomial.

%
%
\section{Reduction to the tree polynomial}

Suppose that $n_0,n_1,\cdots,n_{2k}$ are positive odd integers.
Then we will define an integer $T_k(n_0,\cdots,n_{2k})$. We will
then show that this integer is given by a homogeneous polynomial
in the variables $n_0,\cdots,n_{2k}$ with nonnegative integer
coefficients. We call this the \emph{tree polynomial}. We will
also give a formula for the numbers $b_\ll^n$ in terms of these
polynomials.

\begin{defn}
Let $Sh_k(n_0,\cdots,n_{2k})$ be the set of permutations $\s$ of
the numbers $1,2,\cdots,n$, where $n=\sum n_i$, so that
\begin{enumerate}
    \item $\s(1)=1$,
    \item $\s(n_i+1)<\s(n_i+2)<\cdots<\s(n_{i+1})$ for
    $i=-1,\cdots,2k-1$ where $n_{-1}=0$ and
    \item $\s(n_i+1)<\s(j)<\s(n_{i+1})$ only when $j>n_i$.
\end{enumerate}
We call these permutations \emph{cyclic shuffles}.
\end{defn}

Cyclic shuffles can be described as follows. Take the letters
$a_1,a_2,\cdots,a_{n_0}$ in that order. Then insert the letters
$b_1,b_2,\cdots,b_{n_1}$ in one block between two of the $a$'s or
after the last $a$. There are $n_0$ ways to do this. Next, insert
the letters $c_1,c_2,\cdots,c_{n_3}$ in one block between two
letters in the sequence so far or after the last letter. There are
$n_0+n_1$ ways to do this. Thus the number of elements in this set
is
\[
    |Sh_k(n_0,\cdots,n_{2k})|=n_0(n_0+n_1)(n_0+n_1+n_2)\cdots
    (n_0+\cdots+n_{2k-1}).
\]
Cyclic shuffles have several signs associated to them. The
ordinary sign of $\s$ will be called its \emph{orientation}. We
also have the \emph{selected sign} denoted
$\sgn_\s(a_i,b_j,\cdots)$ which are the sign of $\s$ restricted to
a subset given by selecting one letter of each kind. For example,
take the cyclic shuffle
\[
    \s=a_1a_2b_1c_1c_2c_3c_4b_2a_3.
\]
The orientation is $\sgn(\s)=-1$ and there are $3\cdot2\cdot4$
selected signs
\[
    \sgn_\s(a_i,b_j,c_k)=(-1)^{(j=2)},
\]
ie, the selected sign is negative iff $b_2$ is selected.

The sum of all selected signs will be called the \emph{sign sum}
of $\s$. By the \emph{oriented sign sum} we mean the product of
the sign sum with the orientation of $\s$:
\begin{equation}\label{eq:sign sum for a cyclic shuffle}
    \sgn(\s)\sum \sgn_\s(a_i,b_j,\cdots)=\sgn(\s)\sum \sgn_\s(\s(i),\s(n_0+j),\cdots).
\end{equation}
This has $\prod n_i$ terms. (The sum is for $i=1,\cdots,n_0,
j=1,\cdots,n_1$, etc.) It is easy to see that the oriented sign
sum is divisible by $n_{2k}$ since the selected sign
$\sgn_\s(a_i,b_j,\cdots,y_p)$ is independent of the last index
$p$. (Note that the English language has an even number of letters
so the $2k+1^{\rm st}$ letter cannot be $z$.)

\begin{defn}
Let $T_k(n_0,\cdots,n_{2k})$ be the sum over all cyclic shuffles
$\s$ of the oriented sign sum of $\s$:
\[
    T_k(n_0,\cdots,n_{2k})=\sum_{\s\in Sh_k(n_0,\cdots,n_{2k})}
    \sgn(\s)
    \sum \sgn_\s(a_i,b_j,\cdots)
\]
Let
\[
    Q_k(n_0,\cdots,n_{2k})=\frac{T_k(n_0,\cdots,n_{2k})}
    {|Sh_k(n_0,\cdots,n_{2k})|}
\]
be the average (expected value) of the oriented sign sum over all
cyclic shuffles $\s$.
\end{defn}

We call $T_k(n_0,\cdots,n_{2k})$ the \emph{tree polynomial} since
it is a homogeneous polynomial in $n_0,\cdots,n_{2k}$ with
nonnegative integer coefficients. (Theorem \ref{cor:tree poly is
homogeneous} below). In Section~8 we will show that this
polynomial is in fact the generating function for a statistic on
the set of increasing trees with labels $0,\ldots,2k$. First we
record some obvious properties of the tree polynomial.

\begin{prop}\label{prop:obvious properties of tree polynomial}
For all positive odd integers $n_0,\cdots,n_{2k}$ we have:
\begin{enumerate}
    \item $T_k(n_0,\cdots,n_{2k})$ is an integer.
    \item $T_k(n_0,\cdots,n_{2k})=n_{2k}T_k(n_0,\cdots,n_{2k-1},1)$.
    \item $T_k(n_0,\cdots,n_{2k})$ is divisible by $n_0$ and the
    quotient $T_k(n_0,\cdots,n_{2k})/n_0$ is the sum of $\sgn(\s)
    \sum \sgn_\s(a_i,b_j,\cdots)$ over all cyclic shuffles $\s$ which
    insert the $b$'s after the $a$'s.
\end{enumerate}
\end{prop}

This allows us to compute the first nontrivial tree polynomial.
(The trivial case is $T_0(n_0)=Q_0(n_0)=n_0$.)

\begin{cor}\label{cor:T1}
$T_1(n_0,n_1,n_2)=n_0(n_0+n_1)n_2$, so $Q_1(n_0,n_1,n_2)=n_2$.
\end{cor}

\proof
Since $T_1(n_0,n_1,n_2)=n_2 T_1(n_0,n_1,1)$ it suffices to show
that
\[
    \frac{T_1(n_0,n_1,1)}{n_0}=n_0+n_1.
\]
By Proposition \ref{prop:obvious properties of tree
polynomial}(3), this is given by
$$
    \frac{T_1(n_0,n_1,1)}{n_0}=
    \sum_{i=1}^{n_0}(-1)^i(n_0-2i)n_1
    +\sum_{j=1}^{n_1}(-1)^{j+1}n_0(2j-n_1)=n_1+n_0.
\eqno{\qed}$$

The following theorem tells us that the numbers $b_\ll^n$ (and
thus all $b_\ll^\mu$ and $a_\mu^\ll$) are determined by the tree
polynomials.

\begin{thm}\label{thm:recursive formula for blln}
$b_{\ll,k}^n$ is equal to the sum
\[
    b_{\ll,k}^n=\sum_{(m_0,\cdots,m_{2k})}b_\ll^\mu
    \frac{
    (2m_0+1)Q_k(2m_0+3,2m_1+1,\cdots,2m_{2k}+1)
    }{
    (2m_0+3)(-2)^{k+1}(2k-1)!!
    }
\]
where the sum is over all $2k+1$ tuples of nonnegative integers
$(m_0,\cdots,m_{2k})$ which add up to $n-k$ and $\mu$ is the
partition of $n-k$ given by the nonzero $m_i$.
\end{thm}

\begin{eg}\label{eg:formula for bll 1}
When $k=1$ this formula becomes
\begin{equation}\label{eq:recursive formula for bll 1}
    b_{\ll,1}^n=\sum_{\begin{matrix}%
    \scriptstyle a+b+c=n-1\\
    \scriptstyle a,b,c\geq0\end{matrix}}%
    b_\ll^{[a,b,c]}\frac{(2a+1)(2c+1)}{(2a+3)4}
\end{equation}
where $[a,b,c]$ denotes the multiset $\{a,b,c\}$ with the zero's
deleted. For example, if $\ll=1$ there are three terms with
$[a,b,c]=[1,0,0]=\{1\}$ and (\ref{eq:recursive formula for bll 1})
is
\[
    b_{1,1}^2=\frac{b_1^1}{4}\left(%
\frac35+\frac13+\frac33
    \right)=\frac{29}{60}b_1=\frac{29}{720}.
\]
\end{eg}

\begin{proof} The number $b_{\ll,k}^n$ is given by evaluating the
cup product $\ktilde_\ll\cup\ktilde_k$ on a dual cell of any graph
$\Gam_{2n}$ (with natural orientation) in the Kontsevich cycle
$W_{2n}$. This is given by
\[
    \sum
    o_1o_2\ctilde_\ll(\Gam_0,\cdots,\Gam')
    \ctilde_k(\Gam',\cdots,\Gam_{2n})
\]
where $o_1=o(\Gam_0,\cdots,\Gam')$,
$o_2=o(\Gam',\cdots,\Gam_{2n})$ are the orientations of the front
and back face of the $2n$--simplex $(\Gam_0,\cdots,\Gam_{2n})$. The
sum over all sequences $(\Gam_0,\cdots,\Gam')$ times $o_1$ is the
dual cell of $\Gam'$:
\[
    \sum o_1(\Gam_0,\cdots,\Gam')=D(\Gam').
\]
Consequently, \begin{equation}\label{eq:bll mu}
    \sum o_1\ctilde_\ll(\Gam_0,\cdots,\Gam')=b_\ll^\mu
\end{equation}
if $\Gam'$ lies in the Kontsevich cycle $W_\mu$.

For the other factor, we note that the adjusted cyclic set cocycle
$\ctilde_k$ is a sum of two terms, one for each of the two
vertices of $\Gam'=\Gam_{2n-2k}$ which collapse to a point in the
next graph $\Gam_{2n-2k}$. Each of these vertices gives a pointed
$2k$--simplex. For each such pointed $2k$--simplex, let $v_0,v_1$ be
the two vertices which collapse at the first step and let
$v_2,\cdots,v_{2k}$ be the other vertices of $\Gam'$, indexed
according the order in which they merge with $v_0$.

Since $\Gam'$ must lie in a Kontsevich cycle $W_\mu$, its vertices
$v_i$ must have codimensions $2m_i$ with $m_i\geq0$ so that the
nonzero $m_i$ make up the parts of the partition $\mu$. For each
such sequence $(m_0,\cdots,m_{2k})$ we get a subtotal
\begin{multline*}
    \sum o_1\ctilde_k(\Gam',\cdots,\Gam_{2n})=\\
    \left(
    \frac{2n+3}{2m_0+3}
    \right)
    \left(
    \frac{(2m_0+1)T_k(2m_0+3,2m_1+1,\cdots,2m_{2k}+1)}
    {
    (-2)^{k+1}(2k-1)!!(2m_0+3)(2m_0+2m_1+4)\cdots(2n+3)
    }
    \right)\\
=
    \left(
    \frac{2m_0+1}{2m_0+3}
    \right)
    \left(
    \frac{Q_k(2m_0+3,2m_1+1,\cdots,2m_{2k}+1)}
    {
    (-2)^{k+1}(2k-1)!!
    }
    \right)
\end{multline*}
since there is a $(2n+3)$-to-$(2m_0+3)$ correspondence between
pointed $2k$--simplices and cyclic shuffles. Combine this with
(\ref{eq:bll mu}) and sum over all sequences $(m_0,\cdots,m_{2k})$
to get the result.
\end{proof}

Example \ref{eg:formula for bll 1} allows us to obtain a recursive
formula for $b_{1^n}^n$.

\begin{cor}\label{cor:recursive formula for b1n n}
For all positive $n$ we have
\[
    b_{1^n}^n=4^{-n}n!h(n)
\]
where $h(n)$ is given recursively by $h(0)=1$ and
\[
    h(n+1)=\sum_{\begin{matrix}%
   \scriptstyle a+b+c=n\\
   \scriptstyle a,b,c\geq0\end{matrix}}%
    h(a)h(b)h(c)\frac{(2a+1)(2c+1)}{(2a+3)(n+1)}.
\]
\end{cor}

\begin{proof} In the recursion (\ref{eq:recursive formula for bll
1}) we note that, by the sum of products formula for $b_\ll^\mu$,
we have
\[
    b_{1^n}^{[a,b,c]}=\frac{n!}{a!\, b!\, c!}f(a)f(b)f(c)
\]
where $f(n)=b_{1^n}^n$ for $n\geq1$ and $f(0)=1$. Then
(\ref{eq:recursive formula for bll 1}) becomes
\[
    f({n+1})=\sum_{\begin{matrix}%
    \scriptstyle a+b+c=n\\
    \scriptstyle a,b,c\geq0\end{matrix}}%
\frac{n!}{a!\, b!\, c!}f(a)f(b)f(c)
    \frac{(2a+1)(2c+1)}{(2a+3)4}.
\]
Substitute $f(n)=4^{-n}n!h(n)$ to get the recursion for $h(n)$.
\end{proof}

\begin{eg}
\begin{align*}
    h(1)&=\frac13,&&b_1^1=\frac1{12}\\
    h(2)&=\frac{29}{90},&&b_{11}^2=\frac{29}{720}\\
    h(3)&=\frac{263}{630},&&b_{111}^3=\frac{263}{6720}\\
    h(4)&=\frac{23479}{37800},&&b_{1111}^4=\frac{23479}{403200}
\end{align*}
\end{eg}

The value of $b_{111}^3$ allows us to compute the expansion of
$[W_{111}^3]$ as conjectured by Arbarello and Cornalba
\cite{[Arbarello-Cornalba:96]} and promised in
\cite{[I:GraphCoh]}.

\begin{cor}\label{cor:W111}
$[W_{111}^\ast]=288\ktilde_1^3+4176\ktilde_2\ktilde_1+20736\ktilde_3$
\end{cor}

\proof
By the sum of products formula we have
\[
    b_{111}^{21}=3b_{11}^2b_1^1=3\cdot\frac{29}{720}\cdot
    \frac1{12}=\frac{29}{2880}
\]
\[
    b_{21}^{21}=b_2b_1=\frac1{-120\cdot12}=-\frac1{1440}.
\]
By Equation (\ref{eq:an1 and bn1}) in the introduction which was
proved in \cite{[I:GraphCoh]} but which also follows from Example
\ref{eg:formula for bll 1} above, we have
\[
    b_{21}^3=-\frac{19}{3360}.
\]
Therefore, the coefficients of the expansion
\[
[W_{111}^\ast]=a_{111}^{111}\ktilde_1^3+a_{111}^{21}\ktilde_2\ktilde_1
+a_{111}^3\ktilde_3
\]
are given by
\[
    a_{111}^{111}=\frac{12^3}{3!}=288
\]
\[
    a_{111}^{21}=-\frac{a_{111}^{111}b_{111}^{21}}{b_{21}^{21}}=4176
\]
$$
    a_{111}^3=
    -\frac{a_{111}^{21}b_{21}^3+a_{111}^{111}b_{111}^3}{b_3}=20736.
\eqno{\qed}$$

%
%

\section{First formula for $T_k$}

We will compute the tree polynomial in the case when most of the
entries are equal to 1.

\begin{thm}\label{thm:Tk(n,1,1,...,1)}
\begin{equation*}
    T_k(n, 1, ... , 1,m)=(2k-1)!! mn(n+1)(n+3)...(n+2k-1)
\end{equation*}
\end{thm}

\begin{proof}
Dividing by $n(n+1)(n+2)\cdots(n+2k-1)$ and restricting to the
case $m=1$, it suffices to show that
\begin{equation}\label{eq:Qk(n,1,...,1)}
    {Q_k(n, 1, ... , 1,1)}
    = \frac{(2k-1)!!n!!}{(n+2k-2)!!}
\end{equation}
But $Q_k(n,1,\cdots,1)$ is the expected value of the oriented sign
sum
\[
    \sgn(\s)\sum \sgn_\s(a_i,b_1,b_2,\cdots,b_{2k})
\]
for a random cyclic shuffle $\s$. Since any change in the order of
the $b_i$ leaves this sum invariant, we may assume that the $b$'s
are in correct cyclic order. By Proposition \ref{prop:obvious
properties of tree polynomial}(3), we may also assume that
$b_{2k}$ comes after all the $a$'s. Cyclic shuffles of this kind
are in 1--1 correspondence with ordinary shuffles of
$a_1,\cdots,a_n$ with $b_1,\cdots,b_{2k-1}$ whose oriented sign
sums have expectation values tabulated in the lemma below:
\[
    Q_k(2j-1,1,\cdots,1)=E_0(2j-1,2k-1)=
    \frac{(2j-1)!!(2k-1)!!}{(2j+2k-3)!!}
\]
This gives (\ref{eq:Qk(n,1,...,1)}) proving the theorem.
\end{proof}

\begin{lem}\label{lem:table of E(oriented sign sums)}
Consider all shuffles $\s$ of $a_1,\cdots,a_n$ with
$b_1,\cdots,b_m$ where $n,m$ are nonnegative integers. Then the
sum of the oriented sign sum
\[
    X_0(n,m)=\sum_{\s}\sgn(\s)\sum_{i=1}^n\sgn_\s(a_i,b_1,\cdots,b_m)
\]
and its expected value
\[
    E_0(n,m)=\frac{X_0(n,m)}{\binom{n+m}{n}}
\]
depend on the parity of $n,m$ and are given in the following
table.
\begin{equation*}
\begin{matrix}
 n& m & X_0(n,m) &E_0(n,m)\\
\hline\\
2j & 2k & 2j\binom{j+k}{j}& {\frac{(2j)(2j-1)!!(2k-1)!!}{(2j+2k-1)!!}}\\
\ \\
2j & 2k-1 & 0 & 0\\
\ \\
2j-1 & 2k & {(2j+2k-1)}\binom{j+k-1}{k}
& \frac{(2j-1)!!(2k-1)!!}{(2j+2k-3)!!}\\
\ \\
2j-1 & 2k-1 & 2k\binom{j+k-1}{k} &
\frac{(2j-1)!!(2k-1)!!}{(2j+2k-3)!!}
\end{matrix}
\end{equation*}
\end{lem}

\begin{proof}
Shuffles $\s$ are in 1--1 correspondence with the ways of writing
$m$ as the sum of an $n+1$--tuple of nonnegative integers:
\[
    m=m_0+m_1+\cdots+m_n.
\]
(The corresponding shuffle is $b^{m_0}a_1b^{m_1}a_2\cdots
a_nb^{m_n}$.) The terms in the oriented sign sum are the product
of
\[
    \sgn(\s)=(-1)^{m_{n-1}+m_{n-3}+m_{n-5}+\cdots}
\]
\[
    \sgn_\s(a_i,b_1,b_2,\cdots,b_m)=(-1)^{m_0+m_1+\cdots+m_{i-1}}.
\]
Note that there are $j=\roof{\frac{n}2}$ terms $m_i$ in the
exponent for $\sgn(\s)$. And, when $i$ is even,
$\sgn(\s)\sgn_\s(a_i,b_1,\cdots,b_m)$ has the same form. Thus
\[
    E(\sgn(\s))=E(\sgn(\s)\sgn_\s(a_{2i},b_1,\cdots,b_m)).
\]
Similarly, $\sgn(\s)\sgn_\s(a_{odd},b_1,\cdots,b_m)$ is equal to
$-1$ to the power a sum of $j+(-1)^n$ terms $m_i$ so its expected
value is independent of the subscript of $a$ which can have
$\roof{\frac{n}2}$ different values. The expected value of the
oriented sign sum is thus given by:
\[
    E_0(n,m)=\floor{\frac{n}2}E(\sgn(\s))+
    \roof{\frac{n}2}E(\sgn(\s)\sgn_\s(a_{1},b_1,\cdots,b_m))
\]
\[
    =\frac1{\binom{n+m}{n}}
    \left(\floor{\frac{n}2}\sum_\s\sgn(\s)+
    \roof{\frac{n}2}\sum_\s\sgn(\s)\sgn_\s(a_{1},b_1,\cdots,b_m)
    \right).
\]
The lemma now follows from the following eight computations.
\begin{equation}\label{eq:chart A}
\begin{matrix}
 n& m & \sum\sgn(\s)&\sum\sgn(\s)\sgn_\s(a_{1},b_1,\cdots,b_m)\\
\hline\\
2j & 2k & \binom{j+k}{j}& \binom{j+k}{j}\\
\ \\
2j & 2k-1 & \binom{j+k-1}{j} & -\binom{j+k-1}{j}\\
\ \\
2j-1 & 2k & \binom{j+k-1}{k} & 2\binom{j+k}{k}-\binom{j+k-1}{k}\\
\ \\
2j-1 & 2k-1 & 0 & 2\binom{j+k-1}{j}
\end{matrix}
\end{equation}
We verify the entries in this table starting at the bottom left.
When both $n,m$ are odd there is a fixed point free involution on
the set of shuffles given by switching $m_{2i}\leftrightarrow
m_{2i-1}$. This always changes the sign of $\s$ so
$\sum\sgn(\s)=0$.

When $n=2j-1,m=2k$ we have
\[
    \sgn(\s)=(-1)^{m_0+m_2+\cdots+m_{2j-2}}.
\]
Take the involution on the set of shuffles given as follows. Take
the largest $i$ so that $m_{2i}+m_{2i+1}$ is odd and switch
$m_{2i}\leftrightarrow m_{2i+1}$. If these sums are all even then
take the largest $i$ so that $m_{2i}$ is nonzero. If it is even,
subtract 1 from $m_{2i}$ and add 1 to $m_{2i+1}$. If $m_{2i}$ is
odd, add 1 to it and subtract 1 from $m_{2i+1}$. This sign
reversing involution does not contain all shuffles. The remaining
ones have $m_{even}=0$ and $m_{odd}$ are all even. These shuffles
have positive sign and there are $\binom{j+k-1}{k}$ of them.

The other term on the third line is the sum of
\[
    \sgn(\s)\sgn_\s(a_1,b_1,\cdots,b_m)=(-1)^{m_2+m_4+\cdots+m_{2j-2}}
\]
Here we apply the involution above to $m_2,\cdots,m_{2j-1}$. The
remaining terms are all positive and have $m_2=m_4=\cdots=0$ and
$m_3,m_5,\cdots$ all even and $m_0,m_1$ are arbitrary (with even
sum). There are $\binom{j+k}{k}$ such terms where both $m_0,m_1$
are even. If they are unequal we can subtract 1 from the larger
and add 1 to the smaller, making them both odd. This however
overcounts the terms where $m_0,m_1$ are odd and equal. Thus there
are
\[
    \binom{j+k}{k}-\binom{j+k-1}{k}
\]
terms where both $m_0,m_1$ are odd. The term $\binom{j+k-1}{k}$
counts shuffles where $m_0=m_1$ are odd or even, both kinds being
overcounted once.

By symmetry (switch $n\leftrightarrow m$ and $j\leftrightarrow k$)
we get $\sum\sgn(\s)$ for $n=2j,m=2k-1$. Since $m=2k-1$ is odd,
\[
    \sgn(\s)=(-1)^{m_1+m_3+\cdots+m_{2j-1}}
    =-(-1)^{m_0+m_2+\cdots+m_{2j}}.
\]
This accounts for the $-\binom{j+k-1}{j}$ in the chart. The
remaining three terms are similar.
\end{proof}

\begin{lem}\label{lem:table of E0(oriented sign sums)}
Consider all shuffles $\s$ of $a_0,\cdots,a_n$ with
$b_1,\cdots,b_m$ so that $a_0$ stays on the left (ie, $m_0=0$).
Then the sum $X_1(n,m)$ and average $E_1(n,m)$ of the oriented
sign sum
\[
    \sgn(\s)\sum_{i=0}^n\sgn_\s(a_i,b_1,\cdots,b_m)
\]
are given by
\begin{equation*}
\begin{matrix}
 n& m & X_1(n,m) &E_1(n,m)\\
\hline\\
2j & 2k & (2j+1)\binom{j+k}{j}& {\frac{(2j+1)!!(2k-1)!!}{(2j+2k-1)!!}}\\
\ \\
2j & 2k-1 & \binom{j+k-1}{j} & {\frac{(2j-1)!!(2k-1)!!}{(2j+2k-1)!!}}\\
\ \\
2j-1 & 2k & {(2j+2k)}\binom{j+k-1}{k}
& \frac{(2j+2k)(2j-1)!!(2k-1)!!}{(2j+2k-1)!!}\\
\ \\
2j-1 & 2k-1 & 2k\binom{j+k-1}{k} &
\frac{(2j-1)!!(2k-1)!!}{(2j+2k-3)!!}
\end{matrix}
\end{equation*}
\end{lem}

\begin{rem}
Note that Proposition \ref{prop:obvious properties of tree
polynomial} (3) can be rephrased as:
\[
    Q_k(2j+1,1,\cdots,1)=E_1(2j,2k)=E_0(2j+1,2k-1).
\]
\end{rem}

\begin{proof}
The shuffles in this lemma are the same as those in Lemma
\ref{lem:table of E(oriented sign sums)}. The only difference is
that the oriented sign sum has one more term. The extra term is
\[
    \sgn(\s)\sgn_\s(a_0,b_1,\cdots,b_m)=\sgn(\s).
\]
Therefore
\[
    E_1(n,m)=E(n,m)+E(\sgn(\s)).
\]
The first term is given by Lemma \ref{lem:table of E(oriented sign
sums)}. The second term is given by the first column of
(\ref{eq:chart A}) divided by $\binom{n+m}{n}$.
\end{proof}

\begin{lem}\label{lem:table of E2(oriented sign sums)}
Consider all shuffles $\s$ of $a_0,\cdots,a_{n+1}$ with
$b_1,\cdots,b_m$ so that $a_0$ is the first letter and $a_{n+1}$
is the last. Let $X_2(n,m)$ and $E_2(n,m)$ be the sum and average
value of the oriented sign sum
\[
    \sgn(\s)\sum_{i=0}^{n+1}\sgn_\s(a_i,b_1,\cdots,b_m).
\]
Then
\begin{equation*}
\begin{matrix}
 n& m & X_2(n,m) &E_2(n,m)\\
\hline\\
2j & 2k & (2j+2)\binom{j+k}{j}& {\frac{(2j+2)(2j-1)!!(2k-1)!!}{(2j+2k-1)!!}}\\
\ \\
2j & 2k-1 & 0 & 0\\
\ \\
2j-1 & 2k & {(2j+2k+1)}\binom{j+k-1}{k}
& \frac{(2j+2k+1)(2j-1)!!(2k-1)!!}{(2j+2k-1)!!}\\
\ \\
2j-1 & 2k-1 & -2k\binom{j+k-1}{k} &
-\frac{(2j-1)!!(2k-1)!!}{(2j+2k-3)!!}
\end{matrix}
\end{equation*}
\end{lem}

\begin{proof}
The shuffles in this lemma are the same as those in Lemma
\ref{lem:table of E(oriented sign sums)} with $a_0$ added on the
left and $a_{n+1}$ added on the right. The $a_{n+1}$ on the right
changes the sign by $(-1)^m$ and there are two extra terms in the
oriented sign sum given by
\[
    \sgn(\s)\left(%
    \sgn_\s(a_0,b_1,\cdots,b_m)+\sgn_\s(a_{n+1},b_1,\cdots,b_m)
    \right)%
    =(1+(-1)^m)\sgn(\s).
\]
Therefore
\[
    E_0(n,m)=(-1)^m E(n,m)+(1+(-1)^m)E(\sgn(\s)).
\]
The first term is given by Lemma \ref{lem:table of E(oriented sign
sums)}. The second term is given by the first column of
(\ref{eq:chart A}) times $1+(-1)^m$ divided by $\binom{n+m}{n}$.
\end{proof}

\begin{thm}\label{thm:tree polynomial in main example}
If $p+q=2k-1$ and $n=2r+1$ we get:
\begin{multline*}
    T_k(3,1^p,n,1^q)=
\sum_{s=0}^{\lceil q/2\rceil}
{\frac{q!}{(q-2s+1)!}\binom{r-1+s}{s}(2k-2s)!\, 3\,(k-s+1) \,\cdot
} \\
\left[(q-2s+1)(2r+2s+1)-2s(2k-2s+3)\rule{0pt}{10pt}\right]\nonumber
\end{multline*}
where we use the notation
\[
    T_k(3,1^p,n,1^q)=
    T_k(3,\overbrace{1,\ldots,1}^p,n,\overbrace{1,\ldots,1}^q).
\]
\end{thm}

\begin{proof}
Take cyclic shuffles $\s$ of
\[
    a_1 a_2 a_3 b^1 b^2\cdots b^p c_1\cdots c_n d^1\cdots d^q
\]
Then
\[
    T_k(3,1^p,n,1^q)=\sum_\s\sgn(\s)\sum_{i=1}^3\sum_{j=1}^n
    \sgn_\s(a_ib^1\cdots b^p c_j d^1\cdots d^q)
\]
The cyclic shuffle $\s$ permutes the $b$'s and shuffles them with
the $a$'s, inserts $c_1\cdots c_n$ as one block, then permutes the
$d$'s and shuffles them in.

However, any permutation of the $b$'s will not change the oriented
sign sum since it changes both the orientation and the sign sum by
the same sign. Therefore, it suffices to consider those $\s$ which
do not permute the $b's$ and multiply the result by $p!$.

Similarly, we may assume that the $q$ are in a fixed order, so
that $q^1,\cdots,q^\l$ are shuffled between the $c$'s and
$d^{\l+1},\cdots,d^q$ are shuffled into the $a$'s and $b$'s. The
resulting sum should be multiplied by $q!$. The oriented sign sum
for the shuffles of the $d$'s between the $c$'s is $X_2(n-2,\l)$.
The value of this terms and the remaining factrors depends on the
parity of $\l$.

\textbf{Case 1}\qua $\l=2s$ is even\qua In this case there are an odd number
of letters and an odd number of kinds of letters in the set
$$S=\{c_1,\cdots,c_n,d^1,\cdots,d^\l\}.$$ So, $S$ behaves like a
single letter and we have $p+1+q-\l=2k-2s$ letters shuffled
together in $(2k-2s)!/(p!(q-\l)!)$ ways and then shuffled with
$a_1,a_2,a_3$ keeping $a_1$ first. The contribution to the tree
polynomial given by these shuffles is then
\begin{equation}\label{eq:tree plynomial, even el}
    p!q!
    \frac{(2k-2s)!}{p!(q-\l)!}
    X_1(2,2k-2s)X_2(2r-1,2s).
\end{equation}
By Lemmas \ref{lem:table of E(oriented sign sums)} and
\ref{lem:table of E2(oriented sign sums)} this is equal to
\begin{equation}\label{eq:tree poly, even el expanded}
    \frac{q!(2k-2s)!}{(q-2s)!}
    3(k-s+1)(2r+2s+1)\binom{r+s-1}{s}
\end{equation}

\textbf{Case 2}\qua $\l=2s-1$ is odd\qua In this case the set $S$ has an
even number of letters and an even number of kinds of letters.
Therefore, $S$ can be placed anywhere with the same effect. Since
there are $3+p+q-\l=2k-2s+3$ remaining letters we multiply by this
factor. There are $p$ $b$'s and $q-\l$ $d$'s shuffled together in
\[
    \binom{p+q-\l}{p}=\frac{(2k-2s)!}{p!(q-\l)!}
\]
ways. So the contribution to the tree polynomial of these shuffles
is
\begin{equation}\label{eq:tree plynomial, odd el}
    p!q!(2k-2s+3)
    \frac{(2k-2s)!}{p!(q-\l)!}
    X_1(2,2k-2s)X_2(2r-1,2s-1).
\end{equation}
By Lemmas \ref{lem:table of E(oriented sign sums)} and
\ref{lem:table of E2(oriented sign sums)} this is equal to
\begin{equation}\label{eq:tree poly, odd el expanded}
   - (2k-2s+3)\frac{q!(2k-2s)!}{(q-2s+1)!}
    3(k-s+1)2s\binom{r+s-1}{s}
\end{equation}
Therefore, $T_k(3,1^p,2r+1,1^q)$ is equal to the sum of
(\ref{eq:tree poly, even el expanded}) for $s=0\cdots
\floor{\frac{q}2}$ and (\ref{eq:tree poly, odd el expanded}) for
$s=1\cdots\roof{\frac{q}2}$. Since (\ref{eq:tree poly, even el
expanded}) and (\ref{eq:tree poly, odd el expanded}) are so
similar we can simplify the sum by adding them together to get
\begin{eqnarray}
\label{eq_T} \sum_{s=0}^{\lceil q/2\rceil}
\lefteqn{\frac{q!}{(q-2s+1)!}\binom{r-1+s}{s}(2k-2s)!\, 3\,(k-s+1)
\,\cdot } \\ &&\hspace{1cm}
\left[(q-2s+1)(2r+2s+1)-2s(2k-2s+3)\rule{0pt}{10pt}\right].\nonumber
\end{eqnarray}
The polynomial on the second line consists of the places where the
$\ell=2s$ and $\ell=2s-1$ terms differ.  The sum now runs from
$s=0$ to $s=\lceil q/2\rceil$, which means we have introduced the
extra terms corresponding to $\ell=-1$ and, when $q$ is odd,
$\ell=q+1$, but both of these are zero.
\end{proof}

%
%

\section{A double sum}

Using the formula (\ref{eq_T}) for the tree polynomial
$T_k(3,1^p,2r+1,1^q)$ we are now in a position to compute the
coefficient $b_{r,k}^{r+k}$ for any $r,k\geq0$. 
This section benefitted greatly from the advice of Christian 
Krattenthaler, who pointed out that the techniques of summation in
an earlier version were unnecessarily complicated.

By Theorem
\ref{thm:recursive formula for blln} we have:
\begin{eqnarray}
    b_{r,k}^{r+k}=\sum b_r\frac{
    (2m_0+1)Q_k(2m_0+3,2m_1+1,\cdots,2m_{2k}+1)
    }{
    (2m_0+3)(-2)^{k+1}(2k-1)!!
    }\nonumber
\\
    =\frac{
(2r+1)Q_k(2r+3,1,\cdots,1)
    }{
(2r+3)a_r(-2)^{k+1}(2k-1)!!
    }+\sum_{p+q=2k-1}\frac{
Q_k(3,1^p,2r+1,1^q)
    }{
3a_r(-2)^{k+1}(2k-1)!!
    }\label{eq:brn}
\end{eqnarray}
where $a_r=1/b_r=(-2)^{r+1}(2r+1)!!$. By Theorem
\ref{thm:Tk(n,1,1,...,1)} we know that
\[
    Q_k(2r+1,1,\cdots,1)=\frac{(2k-1)!!(2r+3)!!}{(2r+2k+1)!!}
\]
so the first term in (\ref{eq:brn}) above is equal to
\[
    \frac{2r+1}{a_r(-2)^{k+1}(2k-1)!!}\frac{(2k-1)!!(2r+1)!!}{(2r+2k+1)!!}=
    \frac{2r+1}{(-2)^{r+k+2}(2r+2k+1)!!}.
\]
\begin{lem}\label{lem:double sum}
$$
    \sum_{p+q=2k-1}
Q_k(3,1^p,2r+1,1^q)=
3\frac{(2k+2r+3)}{2k+1}-3\frac{(2r+3)!!\,(2k-1)!!}{(2k+2r+1)!!}
$$
\end{lem}

Suppose for a moment that this is true. Then the second term of
(\ref{eq:brn}) is equal to
\[
    \sum_{p+q=2k-1} \frac{Q_k(3,1^p,2r+1,1^q)}{3a_r(-2)^{k+1}(2k-1)!!}
    =\frac{2r+2k+3}{a_ra_k}-\frac{2r+3}{(-2)^{r+k+2}(2r+2k+1)!!}.
\]
Putting these together we get
\[
    b_{r,k}^{r+k}=\frac{2r+2k+3}{a_ra_k}
    +\frac1{a_{r+k}}
\]
which can be simplified to
\[
    b_{r,k}^{r+k}=b_rb_k(2r+2k+3)+b_{r+k}.
\]
In terms of the adjusted Miller--Morita--Mumford classes this says
\begin{equation}\label{eq:the conjecture MMM in terms of W}
    \ktilde_r\ktilde_k=(b_rb_k(2r+2k+3)+b_{r+k})[W_{r+k}^\ast]
    +\Sym(r,k)b_rb_k[W_{r,k}^\ast]
\end{equation}
where $\Sym(r,k)$ is 2 for $r=k$ and $1$ for $r\noteq k$. Solving
the equation
\[
    a_{r,k}^{r+k}b_{r+k}+a_{r,k}^{r,k}b_{r,k}^{r+k}=0
\]
in which $a_{r,k}^{r,k}=a_ra_k/\Sym(r,k)$ we see that the inverse
coefficient $a_{r,k}^{r+k}$ is given by:%
\begin{multline*}
    a_{r,k}^{r+k}=%
    -\frac{a_ra_k(b_rb_k(2r+2k+3)+b_{r+k})}{\Sym(r,k)b_{r+k}}=
    -\frac{a_ra_k+(2r+2k+3)a_{r+k}}{\Sym(r,k)}\\%
    =
    \frac{(-2)^{r+k+1}}{\Sym(r,k)}(2(2r+1)!!(2k+1)!!-(2r+2k+3)!!)
\end{multline*}
The Kontsevich cycle $W_{r,k}^\ast$ is related to the adjusted MMM
classes by the formula
\[
    [W_{r,k}^\ast]=a_{r,k}^{r+k}\ktilde_{r+k}+a_{r,k}^{r,k}\ktilde_r\ktilde_k.
\]
This gives the following equation as conjectured in
\cite{[I:GraphCoh]}.

\begin{thm}\label{thm:the conjecture Wrk in terms of MMM}
\[
    [W_{r,k}^\ast]=\frac{(-2)^{r+k+1}}{\Sym(r,k)}
    (2(2r+1)!!(2k+1)!!(\ktilde_{r+k}-\ktilde_r\ktilde_k)-(2r+2k+3)!!\ktilde_{r+k})
\]
\end{thm}

\proof[Proof of Lemma \ref{lem:double sum}]
It remains to calculate the sum
$$
\sum_{p+q=2k-1} Q_k(3,1^p,2r+1,1^q),
$$
where in this case
\begin{equation}
\label{eq_Q} Q_k(3,1^p,2r+1,1^q) = \frac{2 (p+2r+3)!}{(p+3)!
(p+q+2r+3)!} T_k(3,1^p,2r+1,1^q).
\end{equation}
By Theorem \ref{thm:tree polynomial in main example} the tree
polynomial $T_k$ is given by
\begin{align} \label{eq_T again}
T_k(3,&1^p,2r+1,1^q)=\nonumber \\
 &\sum_{s=0}^{\lceil q/2\rceil}
\lefteqn{\frac{q!}{(q-2s+1)!}\binom{r-1+s}{s}(2k-2s)!\, 3\,(k-s+1)
\,\cdot } \\ &\hspace{2cm}
\left[(q-2s+1)(2r+2s+1)-2s(2k-2s+3)\rule{0pt}{10pt}\right].\nonumber
\end{align}
We combine equations~(\ref{eq_Q}) and~(\ref{eq_T again}),
eliminating the variable $p=2k-1-q$ and expressing everything in
terms of factorials.  We seek the double sum:
\begin{align}
\label{eq_F}
\sum_{q=0}^{2k-1}& \sum_{s=0}^{\lceil q/2 \rceil} \lefteqn{F(k,r,s,q) \mbox{, where}} \nonumber \\
F(k,r,s,q)& \lefteqn{{}= \frac{6 (2k+2r-q+2)! \,q! \, (r+s-1)! \, (2k-2s)!\,(k-s+1)}
                                             {(2k+2r+2)!\,(2k-q+2)!\,(q-2s+1)! \, s! \, (r-1)!}\,\cdot} \\&\hspace{2cm}
\left[(q-2s+1)(2r+2s+1)-(2s)(2k-2s+3)\rule{0pt}{10pt}\right].\nonumber
\end{align}
The summand $F(k,r,s,q)$ is a hypergeometric term in each of its
variables, so sophisticated summation techniques are available;
see~\cite{A=B} for an introduction.  We are grateful to Christian 
Krattenthaler for suggesting the following path.

The summand $F(k,r,s,q)$ is more manageable as a sum over $q$, so we
will switch the order of the double summation.  The result is {\em almost} 
\begin{equation}
\label{eq_almost}
\sum_{s=0}^{k}\,\, \sum_{q=2s-1}^{2k-1} F(k,r,s,q)
\end{equation}
except that this introduces one new term where $s=0$ and $q=-1$.  Here
$F$ would need delicate handling owing to the $q!$ in its numerator.
We proceed by first calculating the sum in~(\ref{eq_almost}) formally,
and then dealing with the error term that arises when $s=0$.

The inner summation is now over $q$ with $s$ fixed, and the summand is
a $q$--free part times an expression of the form
$$
G(q) = \frac{(A+B-q)!\,q!}{(A-q)!\,(q-C)!}  
  \left[(q-C)(B+C+2)-(C+1)(A-C)\rule{0pt}{10pt}\right],
$$
where $A=2k+2$, $B=2r$, $C=2s-1$.  Gosper's summation algorithm
quickly points out that $G$ has a discrete antiderivative: $G(q) =
H(q+1)-H(q)$, where 
$$
H(q) = \frac{-(A+B+1-q)!\,q!}{(A-q)!\,(q-C-1)\!},
$$
a relation easily verified by hand.  Any definite sum is now easily
computed, and in particular we want $\sum_{q=C}^{A-3} G(q) =
H(A-2)-H(C)$.  Now note that $H(C)=0$ owing to the $(-1)!$ in the
denominator, and we find that the sum is just $$H(A-2) =
\frac{-(A-2)!\,(B+3)!}{2 (A-C-3)!}.$$  Returning to the original
variables and replacing the $q$--free coefficient, we have found that:
\begin{equation}
\label{eq_onedown}
\sum_{q=2s-1}^{2k-1} F(k,r,s,q) =
\frac{-3\, (2k)! \, (2r+3)!}{(r-1)!\,(2k+2r+2)!} \cdot 
  \frac{ (k+1-s) \, (r-1+s)! }{s!}
\end{equation}
Gosper's algorithm reveals that this too has a discrete antiderivative
with respect to $s$:
\begin{eqnarray*}
H(s) &=& \frac{-3\, (2k)! \, (2r+3)!}{(r-1)!\,(2k+2r+2)!} \cdot
\frac{ (kr-rs+2r+k+1)\,(r+s-1)! }{r(r+1)\,(s-1)!} \\[4pt]
&=& \frac{-3\, (2k)! \, (2r+3)! \, (kr-rs+2r+k+1)\,(r+s-1)! }
    {(r+1)!\,(2k+2r+2)!\,(s-1)!}
\end{eqnarray*}
The sum from $s=0$ to $s=k$ is then $H(k+1)-H(0)$, and again we find
that $H(0)=0$.  The full sum is therefore $H(k+1)$, and pulling
factorials together, we get:
\begin{eqnarray}
\label{eq_formalsum}
\sum_{s=0}^{k}\,\, \sum_{q=2s-1}^{2k-1} F(k,r,s,q) &=&
\frac{-6 \, (2k-1)! \, (2r+3)! \, (k+r+1)! }
     {      (k-1)! \, (r+1)! \, (2k+2r+2)! } \nonumber \\
&=&
\frac{-3 \, (2k-1)!! \, (2r+3)!! }{(2k+2r+1)!!}
\end{eqnarray}

Now we compute the error term.  At $s=0$, the left-hand side 
of~(\ref{eq_onedown}) is problematic, but the right-hand side
which we used for further computations is
$$
\frac{-3\, (2k)! \, (2r+3)! \, (k+1)}{(2k+2r+2)!}.
$$
The actual desired value can be easily computed since $F(k,r,0,q)$ is
a $q$--independent factor times the binomial coefficient
$\binom{2k+2r+2-q}{2r}$.
$$
\sum_{q=0}^{2k-1} F(k,r,0,q) =
\frac{6\,(k+1)\,(2k)!}{(2k+2r+2)!} 
  \left( \frac{(2k+2r+3)!}{(2k+2)!} - \frac{(2r+3)!}{2!} \right)
$$
We subtract to find that the error introduced by using the formal
answer at $s=0$ was
\begin{equation}
\label{eq_errorterm}
-\frac{3\,(3+2k+2r)}{(2k+1)}.
\end{equation}
The final answer is the formal sum~(\ref{eq_formalsum}) minus
the error term~(\ref{eq_errorterm}), which is precisely the
conjectured value:
$$
3\frac{(2k+2r+3)}{2k+1}-3\frac{(2r+3)!!\,(2k-1)!!}{(2k+2r+1)!!}
\eqno{\qed}$$

%
%

\section{Reduced tree polynomial}

The second formula for the tree polynomial is based on the
following lemma.

\begin{lem}\label{lem:second counting lemma}
The sum over all sequences of positive integers
\[
    1\leq z_1,z_2,\cdots,z_s\leq n
\]
of the quantity
\[
    (-1)^{z_1+\cdots+z_s}(B(z)-A(z))
\]
where $A(z)$ is the number of positive integers $j$ which are
$\leq z_i$ for an odd number of indices $i$ and $B(z)$ is the
number of positive integers $j\leq n$ so that $j\leq z_i$ for an
even number of $i$ is equal to
\begin{enumerate}
    \item $1$ if $s,n$ are both odd,
    \item $n$ if $s\geq0$ is even and $n$ is odd, and
    \item $\frac{n}{2}(-2)^s$ if $n\geq0$ is even.
\end{enumerate}
\end{lem}

\begin{proof}
First note that this sum can be written as
\begin{equation}\label{eq:for 2nd counting lemma}
    \sum_z(-1)^{\sum z_i}(B(z)-A(z))=\sum_z(-1)^{\sum z_i}\sum_{p=0}^n(-1)^p|L_p(z)|
\end{equation}
where $L_p(z)$ is the set of all $j\in\{1,2,\cdots,n\}$ so that
$j\leq z_i$ for exactly $p$ values of $i$, ie, so that
$z^{-1}[j,n]$ has $p$ elements where $[j,n]$ denotes the set of
integers from $j$ through $n$. This can also be written as
\[
    \sum_{j=1}^n\sum_{z}(-1)^{\sum
    z_i}(-1)^{|z^{-1}[j,n]|}
\]
\[
    =\sum_{j=1}^n\sum_{p=0}^s\binom{s}{p}(-1)^{pj}N(j)^{s-p}N(n-j)^p
\]
where
\[
    N(j)=\sum_{i=1}^{j-1}(-1)^i=\sum_{i=1}^{j+1}(-1)^i=\left\{%
\begin{array}{ll}
    -1 & \hbox{$j$ even} \\
    0 & \hbox{$j$ odd.} \\
\end{array}%
\right.
\]

\textbf{Case 1}\qua $n$ is odd\qua Then either $j$ or $n-j$ is odd
for each $j$. So the summand is nonzero only for $p=0$ or $s$:
\[
    \sum_{p=0}^s\binom{s}{p}(-1)^{pj}N(j)^{s-p}N(n-j)^p
    =N(j)^s+(-1)^{kj}N(n-j)^s
    =(-1)^{s(j+1)}
\]
So the sum (\ref{eq:for 2nd counting lemma}) is equal to
\[
    \sum_{j=1}^n(-1)^{s(j+1)}=\left\{%
\begin{array}{ll}
    1 & \hbox{if $s$ is odd} \\
    n & \hbox{if $s$ is even.} \\
\end{array}%
\right.
\]

\textbf{Case 2}\qua $n$ is even\qua In this case it is possible for
both $j$ and $n-j$ to be even. But then we get
\[
    \sum_{p=0}^s\binom{s}{p}(-1)^{pj}(-1)^s=(-1)^s(1+(-1)^j)^s=(-2)^s.
\]
There are $\frac{n}{2}$ such terms and the other terms, where
$j,n-j$ are both odd, are all zero.
\end{proof}

Using this lemma we get another formula for the tree polynomial
showing that the monomials correspond to increasing trees. Recall
that an \emph{increasing tree} $T$ with vertices $0,1,2,\cdots,2k$
is a tree constructed by attaching the vertices in order. In other
words, $0$ is the root and children are always larger than their
parents. (See \cite{StanleyEC1} for more details about increasing
trees.) For each such $T$ take the monomial in the variables
$x_0,x_1,\cdots,x_{2k}$ given as follows.

For each vertex $i=0,1,\cdots,2k$ of $T$ let $n_i$ be the number
of trees in the forest $T-\{i\}$ with an even number of vertices.
Associate to $T$ the \emph{tree monomial}
\[
    x^T=x_0^{n_0}\cdots x_{2k}^{n_{2k}}.
\]
\begin{eg}\label{ex:tree monomials for k=1}
In the simplest case $k=1$ there are are two increasing trees:
$0-1-2$ and $1-0-2$. The corresponding tree monomials are $x_0x_2$
and $x_1x_2$.
\end{eg}

\begin{lem}\label{lem:homogeneity of tree monomials}
Each tree monomial $x^T$ has degree $2k$.
\end{lem}

\begin{proof}
Since $T$ has an odd number of vertices we can orient each edge so
that it points in the direction in which there are an even number
of vertices. Then $n_i$ is the number of outward pointing edges at
vertex $i$. The sum of the $n_i$ must the number of edges which is
$2k$.
\end{proof}

\begin{thm}\label{thm:tree monomials give tree polynomial}
The sum of the tree monomials $x^T$ is related to the tree
polynomial by:
\[
    x_0\sum_T x^T=T_k(x_0,x_1,\cdots,x_{2k}).
\]
\end{thm}

Suppose for a moment that this is true.

\begin{defn}
We will call
\[
    \widetilde{T}_k(x_0,\cdots,x_{2k}):=\sum_T x^T=\frac1{x_0}T_k(x_0,\cdots,x_{2k})
\]
the \emph{reduced tree polynomial}.
\end{defn}

Since increasing trees are in 1--1 correspondence with permutations
of $1,\cdots,2k$ we get the following.

\begin{cor}\label{cor:tree poly is homogeneous}
The tree polynomial $T_k(x_0,\cdots,x_{2k})$ is a homogeneous
polynomial of degree $2k+1$ with nonnegative integer coefficients
which add up to $(2k)!$, ie, $T_k(1,1,\cdots,1)=(2k)!$.
\end{cor}

Calculations of the reduced tree polynomial tell us something
about permutation. For example, we have the following.

\begin{cor}\label{cor:red tree poly as generating function}
In the special case $x_1=x_2=\cdots=x_{2k}=1$ the reduced tree
polynomial is the generating function
\[
     \wt{T}_k(x, 1, 1, ... , 1)=(2k-1)!!
     (x+1)(x+3)...(x+2k-1)=\sum_{i=0}^k p_ix^i
\]
where $p_i$ is the number of permutations of $2k$ with $i$ even
cycles.
\end{cor}

\begin{proof}
For every increasing tree $T$ the coefficient of $x_0$ in the
monomial $x^T$ is equal to the number of even cycles in the
permutation of $2k$ corresponding to $T$.
\end{proof}

By the following proposition the first variable $x_0$ in the
reduced tree polynomial $\wt{T}_k$ is superfluous.

\begin{prop}\label{prop:tree poly depends only on x0+x1}
The reduced tree polynomial $\wt{T}_k(x_0,\cdots,x_{2k})$ is a
polynomial in the variables $x_0+x_1, x_2, x_3, \cdots, x_{2k}$.
In other words,
\[
    \wt{T}_k(x_0,x_1,\cdots,x_{2k})=
    \wt{T}_k(0,x_0+x_1,x_2,\cdots,x_{2k}).
\]
\end{prop}

\begin{rem}
This means that it suffices to compute $\wt{T}_k$ in the case when
$x_0=0$ since we can recover the general polynomial by
substituting $x_0+x_1$ for $x_1$.
\end{rem}

\begin{proof}
Any increasing tree $T$ contains vertices $0,1$ connected by an
edge together with a certain number of trees $a_1,\cdots,a_r$ with
an odd number of vertices and other trees $b_1,\cdots,b_s$ with an
even number of vertices. These trees can be attached to either $0$
or $1$ giving $2^{r+s}$ different increasing trees. Let $S$ be
this set of increasing trees.

Each $b_i$ gives a factor of either $x_0$ or $x_1$ for $x^T$
depending on whether it is attached to $0$ or $1$. Therefore the
$b_i$'s altogether give a factor of
\[
    (x_0+x_1)^s
\]
to the sum of $x^T$ for all increasing trees in $S$.

The number $r$ must be odd in order for the total number of
vertices to be equal to $2k+1$. Exactly half of the time an odd
number of $a_i$ will be attached to $0$ and the other half of the
time an odd number of $a_i$ will be attached to $1$. Consequently,
the $a_j$ give a factor of
\[
    (x_0+x_1)^{r-1}
\]
to the sum of $x^T$ for all increasing trees in $S$. Thus, the sum
of $x^T$ for all $T\in S$ is equal to $(x_0+x_1)^{r+s-1}$ times a
polynomial in the other variables $x_2,\cdots,x_{2k}$.
\end{proof}

\begin{proof}[Proof of Theorem \ref{thm:tree monomials give tree
polynomial}] Suppose that $\s$ is a cyclic shuffle of the letters
\[
    a^0_1,\cdots,a^0_{x_0},a^1_1,\cdots,a^1_{x_1},\cdots,a^{2k}_1,\cdots,a^{2k}_{x_{2k}}.
\]
Then we associate to $\s$ an increasing tree $T(\s)$ as follows.

To each letter $a^i$ we associate the vertex $i$. We start with
$T_0(\s)$ being just the root $0$ which is associated to
$a^0_1\cdots a^0_{x_0}$. We attach to $T_0(\s)$ the vertex $1$
corresponding to $a^1$. This gives $T_1(\s)$. There are two
possibilities for $T_2(\s)$ as in Example \ref{ex:tree monomials
for k=1}. We get $0-1-2$ if $a^2$ is inserted after (on the right
of) an $a^1_i$. We get $1-0-2$ is $a^2$ is inserted after an
$a^0$. Proceeding by induction suppose that we have constructed
the increasing tree $T_n(\s)$ with vertices $0,1,\cdots,n$. Then
$T_{n+1}(\s)$ is obtained from $T_n(\s)$ by attaching the new
vertex $n+1$ to vertex $j$ if $\s$ inserts $a^{n+1}$ after some
$a^j_i$.

Since there are $x_j$ letters $a^j_i$, the number of cyclic
shuffles $\s$ giving the same increasing tree $T$ is equal to
\[
    x_0^{m_0}x_1^{m_1}\cdots x_{2k-1}^{m_{2k-1}}
\]
where $m_j$ is the number of children that vertex $j$ has.

\medskip
\textbf{Claim}\qua The sum
\begin{equation}\label{eq:sum of oriented sign sums in claim}
    \sum_{T(\s)=T}\sgn(\s)\sum_i\sgn_\s(a^0_{i_0}a^1_{i_1}\cdots a^{2k}_{i_{2k}})
\end{equation}
of the oriented sign sum of $\s$ for all $\s$ with $T(\s)=T$ is
equal to the tree monomial $x^T$ times $x_0$.
\medskip

Since the tree polynomial is the sum of (\ref{eq:sum of oriented
sign sums in claim}) over all increasing trees $T$, this claim
will prove the theorem.

To prove the claim we first consider the unique shuffle $\s_0$
with $T(\s_0)=T$ having the property that each letter is inserted
in the last allowed slot (after the last letter corresponding to
its parent in the increasing tree $T$). The oriented sign sum of
$\s_0$ is equal to the product
\begin{equation}\label{eq:j=0 case of induction}
    \sgn(\s_0)\sum_i\sgn_\s(a^0_{i_0}a^1_{i_1}\cdots
    a^{2k}_{i_{2k}})=x_0x_1\cdots x_{2k}
\end{equation}
since every summand is equal to $1$.

The statement (\ref{eq:j=0 case of induction}) is the base case
($j=0$) of the following induction hypothesis:
\begin{equation}\label{eq:induction hypothesis}
    \sum_\s\sgn(\s)\sum_i\sgn_\s(a^0_{i_0}a^1_{i_1}\cdots
    a^{2k}_{i_{2k}})=
    x_0^{n_0}\cdots x_{j-1}^{n_{j-1}} x_j\cdots x_{2k}
\end{equation}
if the sum is taken over all $\s$ so that
\begin{enumerate}
    \item $T(\s)=\s$
    \item for all $i\geq j$ the children of $i$ are inserted in
    the last allowed slot (after the last $a^i$).
\end{enumerate}We recall that $n_i$ is the number of components of
$T(\s)-i$ having an even number over vertices.

Suppose by induction that (\ref{eq:induction hypothesis}) holds
for $j$. To extend it to $j+1$ we need to allow the children of
vertex $j$ to be inserted at any of the $x_j$ allowed points.

Let $b^1,b^2,\cdots,b^r$ be the letters corresponding to the
children of $j$ with an odd number of descendants. Then each $b^i$
has the property that it, together with all its descendants, can
be moved to any other slot without changing the oriented sign sum.
This is because both the shuffle and the permutation of selected
letters changes by an even permutation. Consequently, the sum
(\ref{eq:induction hypothesis}) is multiplied by $x_j^r$ bringing
the value of (\ref{eq:induction hypothesis}) to
\begin{equation}\label{eq:new product}
    x_0^{n_0}\cdots x_{j-1}^{n_{j-1}} x_j^{r+1}\cdots x_{2k}.
\end{equation}
Let $c^1,\cdots,c^s$ be the other children of $j$, the ones with
an even number of descendants. Let $z_1,\cdots,z_s$ denote the
indices of the letter $a^j$ after which these letters are
inserted, eg, $c^1$ is inserted after $a^j_{z_1}$. Take the sum:
\begin{equation}\label{eq:partial oriented sign sum}
    \sum_z\sgn(\s)\sum\sgn_\s(a^0_{i_0}\cdots a^{2k}_{i_{2k}})
\end{equation}
over all $(x_j)^s$ insertion points $z=(z_1,\cdots,z_s)$ for all
of the children $c^i$ together with their descendants. The
question is: How do the terms in this sum compare to the term in
which all of the $z_i$ are maximal (equal to $x_j$)?
\medskip

\textbf{Case 1}\qua $s$ is odd\qua (Then there are an odd number of
vertices in $T$ minus $j$ and its descendants. So $n_j=r$.) In
this case we claim that the sum (\ref{eq:partial oriented sign
sum}) is equal to $\frac1{x_j}$ times the summand in which each
$z_i$ is maximal. This is the first case of Lemma \ref{lem:second
counting lemma}. To see this consider what happens when we
decrease by one the insertion point $z_i$ of $c^i$ and its
descendants. This will change the sign of $\s$ by $(-1)^{m+1}$
where $m$ is the number of other $b_p$ which are transposed with
$c^i$. But the selected sign $\sgn_\s(a^0_{i_0}\cdots
a^{2k}_{i_{2k}})$ also changes by $(-1)^m$ so the net effect is to
change the sign of the oriented sign sum. Since $x_j$ is odd and
the sign changes $x_j-z_i$ times, this gives a sign factor of
\begin{equation}\label{eq:sign factor in proof of thm}
    (-1)^{s+z_1+\cdots+z_s}=-(-1)^{z_1+\cdots+z_s}.
\end{equation}
For each value of the index $i_j$ of $a^j$, the selected sign only
changes when some $z_i$ goes below $i_j$. Taking the sum over all
values of $i_j$ we get a factor of
\[
    A(z)-B(z)
\]
instead of $x_j$ where $A(z),B(z)$ are as defined in Lemma
\ref{lem:second counting lemma}. This factor, together with
(\ref{eq:sign factor in proof of thm}), adds up to $1$ by the
lemma. This is instead of the factor of $x_j$ which we get in the
case when each $z_i$ is maximal. So, the sum (\ref{eq:induction
hypothesis}) for $j+1$ is equal to
\[
    x_0^{n_0}\cdots x_{j-1}^{n_{j-1}} x_j^{r}\cdots x_{2k}
\]
which is correct since $r=n_j$.
\medskip

\textbf{Case 2}\qua $s$ is even\qua (Then there are an even number of
vertices in $T$ minus $j$ and its descendants making $n_j=r+1$.)
In this case we claim that the sum (\ref{eq:partial oriented sign
sum}) is equal to the term in which all the $z_i$ are maximal. The
proof is the same as in Case 1, using the second case of Lemma
\ref{lem:second counting lemma}. This leaves the product
(\ref{eq:new product}) unchanged. But this is correct since
$n_j=r+1$.\end{proof}

%
%

\section{Recursion for $\wt{T}_k$}

We will show that the reduced tree polynomial $\wt{T}_k$ (in
variables $x_0, \cdots , x_{2k}$) satisfies a recursion which we
can express in terms of an exponential generating function. First
we need to generalize the reduced tree polynomial.

\begin{defn} For $k,n\geq0$ let $L_k^n$ be the polynomial in
generators $x_0,\cdots,x_{2k}$ given by
\[
    L_k^n=\sum_T \frac{x^T}{x_{2k+1}\cdots x_{2k+2n}}
\]
where the sum is taken over all increasing trees $T$ with vertices
$0$ through $2k+2n$ of which the last $2n$ are leaves. To simplify
notation we write the summand above as $\hat{x}^T$ (ie, this is
$x^T$ with $x_{2k+1},\cdots,x_{2k+2n}$ set equal to zero. If we
delete the last $2n$ vertices from $T$ we get what we call the
\emph{base tree} $T_0$ which is an arbitrary increasing tree with
vertices $0,\cdots,2k$.
\end{defn}

We make some trivial observations about this polynomial.

\begin{prop}\label{prop:observations about Lkn}\begin{enumerate}
    \item $L_k^0=\wt{T}_k$ is the reduced tree polynomial.
    \item $L_0^n(x_0)=1$
    \item The polynomial $L_k^n$ has nonnegative integer coefficients adding
up to $$L_k^n(1,\cdots,1)=(2k)!(2k+1)^{2n}.$$
\end{enumerate}
\end{prop}

Let $g_k(t)$ be the exponential generating function:
\[
    g_k(t)=\sum_{n=0}^\infty L_k^n\frac{t^{2n}}{(2n)!}
\]
Then $g_k(0)=L_k^0=\wt{T}_k$. So it suffices to compute $g_k(t)$
for all $k$. When $k=0$ we have $L_0^n=1$ so
\[
    g_0(t)=\sum \frac{t^{2n}}{(2n)!}=\cosh t.
\]

\begin{thm}\label{thm:recursion for gk(t)}
The generating function $g_k(t)$ which gives $g_k(0)=\wt{T}_k$ is
given recursively as follows.
\begin{enumerate}
    \item $g_0(t)=\cosh t$
    \item $g_{k+1}(t)=g_k(t)\left(z_{2k}z_{2k+2}\sinh^2 t%
      +z_{2k} y_2\right) %
     $%
 \[+g_k'(t)z_{2k+1}(y_1+y_2)\sinh t\cosh t+
 g_k''(t)y_1y_2\cosh^2 t\]
\end{enumerate}
where we use the notation: $z_j = x_0+\cdots+x_j, y_i=x_{2k+i}$.
\end{thm}

We will obtain a recursive formula to compute the polynomials
$L_k^n$ and use the recursion to show the theorem. We begin with
the first nontrivial case $k=1$.

For $k=1$ there are two possibilities for the base tree
(consisting of the vertices $0,1,2$). They are connected either as
$1-0-2$ or $0-1-2$. In each case we attach $2n$ leaves in all
$3^{2n}$ possible ways.

Let $\a,\b,\g$ be the number of leaves attached to $1,2,0$
respectively. We note that there are
\[
    \sum_{j=1}^n\binom{2n}{2j}2^{2j-1}=\frac{3^{2n}-1}4
\]
ways for $\a/\b/\g$ to be odd/odd/even and similarly for the cases
odd/even/odd and even/odd/odd. This leaves
    \[
    \frac{3^{2n}+3}4
    \]
ways for $\a,\b,\g$ to be all even. We determine the monomials
$\hat{x}^T$ in each case.
\begin{enumerate}
    \item Base $1-0-2$ with $\a,\b,\g$ all even. In this case the
    monomial is $\hat{x}^T=x_1x_2$. So the contribution is
    \[
  \left( \frac{3^{2n}+3}4\right) x_1 x_2.
    \]
    \item Base $1-0-2$ with $\a,\b$ both odd (and $\g$ even). Then the
    monomial is $\hat{x}^T=x_0^2$. So the contribution is
    \[
   \left(\frac{3^{2n}-1}4\right)x_0^2.
    \]\item Base $1-0-2$ with $\g$ odd. Then the
    monomial is $\hat{x}^T=x_0x_i$ where $i=1,2$ with equal probability. So the contribution is
    \[
   \left(\frac{3^{2n}-1}4\right) x_0(x_1+x_2).
    \]
\end{enumerate}
Adding these three together we get
\begin{equation*}
    \frac{3^{2n}}{4}(x_0(x_0+x_1+x_2)+x_1x_2)+
    \frac14(3x_1x_2-x_0(x_0+x_1+x_2))
\end{equation*}
If the base tree is $0-1-2$ then we just switch $x_0$ and $x_1$ in
the above expression. Adding these two cases gives
\begin{equation}\label{eq:L1n}
    L_1^n=\frac{3^{2n}}4
{(x_0+x_1)(2x_2+x_0+x_1)}
    +
    \frac14{(x_0+x_1)(2x_2-x_0-x_1)}.
\end{equation}
Note that $n$ occurs only in the exponent of $3$. More generally,
we have the following.
\begin{lem}\label{lem:Lkn is lin comb of (2s+1)2n}
\[
    L_k^n=\sum_{s=0}^k4^{-k}(2s+1)^{2n}P_k^{2s+1}
\]
where $P_k^{2s+1}$ is a polynomial in $x_0,\cdots,x_{2k}$ with
integer coefficients depending only on $k,s$.
\end{lem}

\begin{rem}\label{rem:gk(t) in terms of Pks}
This lemma can be rephrased in terms of the exponential generating
function $g_k(t)$ as follows.
\[
    g_k(t)=\sum_{n\geq0}L_k^n\frac{t^{2n}}{(2n)!}=
    \sum_{n,c} \frac{P_k^c}{4^{k}}\frac{c^{2n}t^{2n}}{(2n)!}=
    \sum_{s=0}^k 4^{-k}P_k^{2s+1}\cosh((2s+1)t).
\]
\end{rem}

We will prove Lemma~\ref{rem:gk(t) in terms of Pks} and find a
recursion for $L_k^n$ at the same time. Suppose we know the
polynomial $L_k^n$ for all $n$ and we wish to compute $L_{k+1}^n$.
This is a sum of monomials $\hat{x}^T$. There are again two cases
for the base tree $T_0$. Either $2k+1,2k+2$ are leaves of the base
tree or $2k+2$ is attached to $2k+1$. In both cases we attach $2n$
leaves to $T_0$, $\a$ on $2k+1$, $\b$ on $2k+2$ and $\g$ on $T_-$
where $T_-$ is $T_0$ with the vertices $2k+1,2k+2$ removed.
\medskip

\textbf{Case 1}\qua $2k+1,2k+2$ are leaves of the base tree $T_0$.
\begin{enumerate}
    \item $\a,\b,\g$ all even with $\g=2m$. In this case the
    vertices $2k+1,2k+2$ act like leaves and $T$ looks like $T_-$
    with $2m+2$ leaves. The monomials in this case add up to
    $$L_k^{m+1}
x_{2k+1}x_{2k+2}.$$ We need to multiply this with the number of
choices for the $\a,\b,\g$ leaves which is
    \[
    \binom{2n}{2m}2^{2n-2m-1}
    \]
    if $0\leq m< n$ and $1$ if $m=n$. This gives a contribution
    of
    \begin{equation}\label{eq:first contribution to Lk+1n}
L_k^{n+1} x_{2k+1}x_{2k+2}+
\sum_{m=0}^{n-1}L_k^{m+1}\binom{2n}{2m}2^{2n-2m-1}
x_{2k+1}x_{2k+2}.
    \end{equation}
    \item $\a,\b$ both odd with $\g=2m$. In this case the vertices
$2k+1,2k+2$ simply add a factor of $x_ix_j$ to $\hat{x}^T$ if they
are attached to vertices $i,j\leq2k$. Taking the sum over all
$i,j$ we get a factor of $z_{2k}^2$ where
\[
    z_{2k}=x_0+x_1+\cdots+x_{2k}.
\]
The contribution to $L_{k+1}^n$ is thus
\begin{equation}\label{eq:second contribution to Lk+1n}
    \sum_{m=0}^{n-1}L_k^m\binom{2n}{2m}2^{2n-2m-1}z_{2k}^2.
\end{equation}
    \item $\g=2m-1$. In this case one of $\a,\b$ is odd and the other is
    even. This gives a factor of $z_{2k}(x_{2k+1}+x_{2k+2})$ for a
    contribution of
\begin{equation}\label{eq:third contribution to Lk+1n}
    \sum_{m=1}^{n}L_k^m\binom{2n}{2m-1}2^{2n-2m}z_{2k}(x_{2k+1}+x_{2k+2}).
\end{equation}
\end{enumerate}

\textbf{Case 2}\qua $2k+2$ is attached on $2k+1$.
\begin{enumerate}
    \item $\a,\b,\g$ all even with $\g=2m$. Then the tree
    consisting of vertices $2k+1,2k+2$ and $\a+\b$ leaves has an
    even number of vertices and contributes a factor of
    $z_{2k}x_{2k+2}$. As in Case 1(1) we get a contribution to
    $L_{k+1}^n$ of
    \begin{equation}\label{eq:fourth contribution to Lk+1n}
L_k^{n} z_{2k}x_{2k+2}+
\sum_{m=0}^{n-1}L_k^m\binom{2n}{2m}2^{2n-2m-1} z_{2k}x_{2k+2}.
\end{equation}
    \item $\a,\b$ both odd with $\g=2m$. This time we get a factor
    of $z_{2k}x_{2k+1}$ so the contribution is
    \begin{equation}\label{eq:fifth contribution to Lk+1n}
\sum_{m=0}^{n-1}L_k^m\binom{2n}{2m}2^{2n-2m-1} z_{2k}x_{2k+1}.
\end{equation}
\item $\g=2m-1$. This is just like Case 1(3). The tree with
vertices $2k+1,2k+2$ and $\a+\b$ leaves acts like one leaf. We get
a factor of $x_{2k+1}^2$ or $x_{2k+1}x_{2k+2}$ depending on
whether $\a$ or $\b$ is even. Thus the contribution is
\begin{equation}\label{eq:sixth contribution to Lk+1n}
    \sum_{m=1}^{n}L_k^m\binom{2n}{2m-1}2^{2n-2m}x_{2k+1}(x_{2k+1}+x_{2k+2}).
\end{equation}
\end{enumerate}

The value of $L_{k+1}^n$ is given by adding these six terms:
\[
    L_{k+1}^n=(\ref{eq:first contribution to Lk+1n})+%
(\ref{eq:second contribution to Lk+1n})+%
(\ref{eq:third contribution to Lk+1n})+%
(\ref{eq:fourth contribution to Lk+1n})+%
(\ref{eq:fifth contribution to Lk+1n})+%
(\ref{eq:sixth contribution to Lk+1n}).
\]
To simplify the computation we need to use Lemma \ref{lem:Lkn is
lin comb of (2s+1)2n} and the following two formulas.
\begin{equation*}
    \sum_{m=0}^{n-1}\binom{2n}{2m}c^{2m}2^{2n-2m}=
    \frac{(c+2)^{2n}+(c-2)^{2n}}2-c^{2n}
\end{equation*}
\begin{equation*}
    \sum_{m=1}^n\binom{2n}{2m-1}c^{2m}2^{2n-2m}=\frac{c}2\left(
    \frac{(c+2)^{2n}-(c-2)^{2n}}2\right)
\end{equation*}

\begin{proof}[Proof of Lemma \ref{lem:Lkn is lin comb of
(2s+1)2n}] We know that the lemma holds for $k=0,1$ so suppose
that $k\geq1$ and the lemma holds for $k$. Substituting the
expression $c^{2m}$ for $L_k^m$ and letting $y_i=x_{2k+i}$ we get
the following.
\[
    \text{expression}(\ref{eq:first contribution to Lk+1n})=
    c^{2n+2}y_1y_2+
    \frac{c^2}2\left(
    \frac{(c+2)^{2n}+(c-2)^{2n}}2-c^{2n}
    \right)y_1y_2
\]
\[
    ={c^2}\left(\frac{
    (c+2)^{2n}+(c-2)^{2n}}4+\frac{c^{2n}}2
    \right)y_1y_2
\]
\[
    \text{expression}(\ref{eq:second contribution to Lk+1n})=
     \left(\frac
    {(c+2)^{2n}+(c-2)^{2n}}4-\frac{c^{2n}}2
    \right)z_{2k}^2
\]
\[
    \text{expression}(\ref{eq:third contribution to Lk+1n})=
{c}\left(\frac
    {(c+2)^{2n}-(c-2)^{2n}}4
\right)z_{2k}(y_1+y_2)
\]
\[
    \text{expression}(\ref{eq:fourth contribution to Lk+1n})=
c^{2n}z_{2k}y_2+\frac{1}2\left(
    \frac{(c+2)^{2n}+(c-2)^{2n}}2-c^{2n}
\right)z_{2k}y_2
\]
\[
    =\left(\frac
    {(c+2)^{2n}+(c-2)^{2n}}4+\frac{c^{2n}}2
\right)z_{2k}y_2
\]
\[
    \text{expression}(\ref{eq:fifth contribution to Lk+1n})=
\left(\frac
    {(c+2)^{2n}+(c-2)^{2n}}4-\frac{c^{2n}}2
\right)z_{2k}y_1
\]
\[
    \text{expression}(\ref{eq:sixth contribution to Lk+1n})=
    {c}\left(\frac
    {(c+2)^{2n}-(c-2)^{2n}}4
\right)y_1(y_1+y_2)
\]
Collect together the terms with $c^{2n}/4,(c\pm2)^{2n}/4$. Then,
for every $c^{2n}$ term which occurs in $L_k^n$ we get the
following three terms in $L_{k+1}^n$.
\begin{multline}\label{eq:Lk+1n part zero}
    \frac{c^{2n}}4\left(
    2c^2y_1y_2-2z_{2k}^2+2z_{2k}y_2-2z_{2k}y_1
    \right)
   = \frac{c^{2n}}4\left(
    2c^2y_1y_2-2z_{2k}(z_{2k+1}-y_2)
    \right)
\end{multline}
\begin{multline}\label{eq:Lk+1n part plus}
    \frac{(c+2)^{2n}}4\left(
c^2y_1y_2+z_{2k}^2+cz_{2k}(y_1+y_2)+z_{2k}y_2+z_{2k}y_1+cy_1(y_1+y_2)
    \right)
\\
    =\frac{(c+2)^{2n}}4
    \left(
    c^2y_1y_2+cz_{2k+1}(y_1+y_2)+z_{2k}z_{2k+2}
    \right)
\end{multline}
\begin{multline}\label{eq:Lk+1n part minus}
    \frac{(c-2)^{2n}}4\left(
c^2y_1y_2+z_{2k}^2-cz_{2k}(y_1+y_2)+z_{2k}y_2+z_{2k}y_1-cy_1(y_1+y_2)
    \right)
\\
    =\frac{(c-2)^{2n}}4
       \left(
    c^2y_1y_2-cz_{2k+1}(y_1+y_2)+z_{2k}z_{2k+2}
    \right)
\end{multline}

If $L_k^n$ is a linear combination of $c^{2n}/4^k$ for
$c=1,3,\cdots,2k+1$ then $L_{k+1}^n$ is a linear combination of
the above three expressions which in turn are linear combinations
of $c^{2n}/4^{k+1}$ for $c=1,3,\cdots,2k+3$. This proves the
lemma.
\end{proof}

If we change the sign of $c$ then (\ref{eq:Lk+1n part plus}),
(\ref{eq:Lk+1n part minus}) are interchanged and (\ref{eq:Lk+1n
part zero}) remains the same. Consequently, these three
expressions directly translate into the following recursion for
the coefficients $P_k^c$.

\begin{thm}\label{thm:recursion for Pks}
$L_k^n=\sum_{s=0}^k 4^{-k}P_k^{2s+1}(2s+1)^{2n}$ where
$P_k^c=P_k^{-c}$ is given for all odd integers $c$ as follows.
\[ P_0^1=P_0^{-1}=1,\
P_0^c=0\ \text{if}\ |c|>1
\]
\begin{align*}
P_{k+1}^c=\ &P_k^c(2c^2 y_1 y_2-2 z_{2k}(z_{2k+1}-y_2))\\
    &+P_k^{c-2}\left(
    (c-2)^2y_1y_2+(c-2)z_{2k+1}(y_1+y_2)+z_{2k}z_{2k+2}
    \right)\\
    &+P_k^{c+2}\left(
    (c+2)^2y_1y_2-(c+2)z_{2k+1}(y_1+y_2)+z_{2k}z_{2k+2}
    \right)
\end{align*}
where $z_j = x_0+\cdots+x_j$ and $y_i=x_{2k+i}$.
\end{thm}

\begin{cor}\label{cor:formula for reduce tree poly}
The reduced tree polynomial is given by
\[
    \wt{T}_k=L_k^0=\sum_{s=0}^k 4^{-k}P_k^{2s+1}.
\]
\end{cor}

\begin{proof}[Proof of Theorem \ref{thm:recursion for gk(t)}]
The exponential generating function for $L_k^n$ is
\[
    g_k(t)=\sum_{n,c}\frac{P_k^c}{4^k}c^{2n}\frac{t^{2n}}{(2n)!}=
    \sum_{c}\frac{P_k^c}{4^k}\cosh ct.
\]
Using the hypertrigonometric identity
\[
    \cosh(ct\pm2t)=\cosh ct\cosh 2t\pm\sinh ct\sinh 2t
\]
we get:
\begin{align*}
g_k\cosh 2t&=\sum_c \frac{P_k^c}{4^k}\frac12(
    \cosh(ct+2t)+\cosh(ct-2t))\\%
g_k'\sinh 2t&=\sum_c \frac{P_k^c}{4^k}\frac{c}2(
    \cosh(ct+2t)-\cosh(ct-2t))\\%
g_k''\cosh 2t&=\sum_c \frac{P_k^c}{4^k}\frac{c^2}2(
    \cosh(ct+2t)+\cosh(ct-2t)).
\end{align*}
So, the recursion in Theorem \ref{thm:recursion for Pks} gives us:
\begin{multline*}
    g_{k+1}=\frac{g_k''}{2}y_1y_2-\frac{g_k}2z_{2k}(z_{2k+1}-y_2)\\%
    +\frac{g_k''}{2}(\cosh 2t)y_1y_2
    +\frac{g_k'}{2}(\sinh 2t)z_{2k+1}(y_1+y_2)
    +\frac{g_k}{2}(\cosh 2t)z_{2k}z_{2k+2}
\end{multline*}
Simplify this to get the theorem.
\end{proof}

%
%

\section{Examples of $\wt{T}_k$}

We will use the following version of Theorem \ref{thm:recursion
for Pks} to compute the reduced tree polynomial $\wt{T}_k$ for
small $k$. By Proposition \ref{prop:tree poly depends only on
x0+x1} it suffices to consider the case when $x_0=0$. We use the
following version of the recurrence.

\begin{align*}
P_{k+1}^c&=
P_k^c(2c^2x_{2k+1}x_{2k+2}-2z_{2k}(z_{2k+1}-x_{2k+2}))\\
&+P_k^{c-2}(z_{2k}+(c-2)x_{2k+1})(z_{2k+1}+(c-1)x_{2k+2})\\
&+P_k^{c+2}(z_{2k}-(c+2)x_{2k+1})(z_{2k+1}-(c+1)x_{2k+2})
\end{align*}
\[
    P_0^1=1\qquad\qquad
\wt{T}_0=P_0^1=1
\]
Since $z_0=x_0=0$ we get:
\begin{align*}
    P_1^3&=P_0^1x_1(x_1+2x_2)=x_1^2+2x_1x_2 &&\\
    P_1^1&=P_0^1(2x_1x_2-x_1(x_1))=-x_1^2+2x_1x_2 &&
\end{align*}
\[
    \wt{T}_1=\frac14(P_1^1+P_1^3)=x_1x_2
\]
When $k=2$ the polynomials $P_k^c$ and $\wt{T}_2$ are still
manageable:
\begin{align*}
    P_2^5&=P_1^3(z_2+3x_3)(z_3+4x_4)\\
    &=x_1(x_1+2x_2) (z_2+3x_3)(z_3+4x_4)\\
    P_2^3&=P_1^3(18x_3x_4-2z_2(z_3-x_4))+
    P_1^1(z_2+x_3)(z_3+2x_4)\\
    &=x_1(x_1+2x_2)(18x_3x_4-2z_2(z_3-x_4))+
    x_1(-x_1+2x_2)(z_2+x_3)(z_3+2x_4)\\
    P_2^1&=P_1^1(2 x_3 x_4- 2 z_2 (z_3 -x_4) +(z_2 - x_3) z_3)+
    P_1^3(z_2-3x_3)(z_3-2x_4)\\
    &=x_1(-x_1+2x_2)(2 x_3 x_4- 2 z_2 (z_3 -x_4) +(z_2 - x_3) z_3)
    +x_1(x_1+2x_2)\\
    &\hspace{9cm}(z_2-3x_3)(z_3-2x_4)
\end{align*}
\[
    \wt{T}_2=\frac1{4^2}(P_2^1+P_2^3+P_2^5)=
    x_1^2  x_2 x_4 + x_1 x_2^2  x_4 + 2 x_1^2  x_3 x_4 + 5
x_1 x_2 x_3 x_4
\]
For $k\geq3$ both $P_k^c$ and $\wt{T}_k$ become more complex
(except for $P_k^{2k+1}$):
\begin{align*}
P_3^7=&\ x_1(x_1+2x_2) (z_2+3x_3)(z_3+4x_4)(z_4+5x_5)(z_5+6x_6)\\
    \wt{T}_3=&\ \frac1{4^3}(P_3^1+P_3^3+P_3^5+P_3^7)\\
    =&\ 8 x_1^2  x_2 x_3 x_4 x_6 + 16 x_1^2  x_2 x_3 x_5 x_6 + x_1^3  x_2 x_4 x_6 + 2
x_1^2 x_2^2  x_4 x_6 + x_1^2  x_2 x_4^2  x_6 \\
+& 23 x_1^2  x_2 x_4 x_5 x_6+ 6 x_1 x_2^2  x_3 x_4 x_6 + 12 x_1
x_2^2  x_3 x_5 x_6 + 2 x_1^3  x_3 x_4 x_6 + 2 x_1^2 x_3^2  x_4 x_6
\\
+& 2 x_1^2  x_3 x_4^2 x_6+ 28 x_1^2  x_3 x_4 x_5 x_6 + 5 x_1 x_2
x_3^2  x_4 x_6 + 10 x_1 x_2 x_3^2  x_5 x_6 + 2 x_1^3  x_2 x_5
x_6\\
+& 4 x_1^2  x_2^2  x_5 x_6+ 6 x_1^3  x_4 x_5 x_6 + 4 x_1^2  x_3^2
x_5 x_6 + 5 x_1 x_2 x_3 x_4^2  x_6 + 61 x_1 x_2 x_3 x_4 x_5 x_6
\\
+& x_1 x_2^3  x_4 x_6 + x_1 x_2^2  x_4^2  x_6 + 2 x_1 x_2^3  x_5
x_6 + 4 x_1^3  x_3 x_5 x_6 + 17 x_1 x_2^2 x_4 x_5 x_6
\end{align*}

The coefficients of $\wt{T}_k$ tell us something about increasing
trees. For example, $61$ (the coefficient of $x_1 x_2 x_3 x_4 x_5
x_6$) is the number of increasing trees in which each node has an
even number of children.

%
%

\subsection*{Summary of algorithm}
\addcontentsline{toc}{subsection}{Summary of algorithm}

First we obtain the reduced tree polynomial by substituting
$x_0+x_1$ for $x_1$. For example $\wt{T}_2$ is given by:
\begin{align*}
    \wt{T}_2(x_0,\cdots,x_4)&=\\
(x_0+&x_1)^2  x_2 x_4 + (x_0+x_1) x_2^2  x_4
     + 2 (x_0+x_1)^2  x_3 x_4 + 5(x_0+x_1) x_2 x_3 x_4
\end{align*}
Next, we need to find $Q_k$ which is given in general by
\[
    Q_k(x_0,\cdots,x_{2k})=\frac{\wt{T}_k(x_0,\cdots,x_{2k})}
    {z_1z_2\cdots z_{2k-1}}.
\]
For $k=2$ this is
\[
Q_k(x_0,\cdots,x_{4})=\frac{(x_0+x_1+x_2+x_3)  x_2 x_4 + 2
(x_0+x_1+x_2) x_3 x_4 + 2 x_2 x_3 x_4}
    {(x_0+x_1+x_2)(x_0+x_1+x_2+x_{3})}.
\]
Take any partition $\mu$ of $m$ with at most $2k+1$ parts. Write
the parts in any order and insert $0$ at the end:
\[
    \mu=(m_0,m_1,\cdots,m_{2k}),\qquad \sum m_i=m.
\]
The simplest example has only one part: $\mu=m0^{2k}$. Let
\[
    R_k(\mu):=\frac{2m_0+1}{2m_0+3}Q_k(2m_0+3,2m_1+1,\cdots,2m_{2k}+1).
\]
Let $S_k(\mu)$ be the symmetrized version of $R_k$:
\[
    S_k(\mu):=\frac1{\Sym(\mu)}\sum_\s
    R_k(m_{\s(0)},m_{\s(1)},\cdots,m_{\s(2k)}),
\]
where the sum is over all permutations $\s$ of the letters
$0,\cdots,2k$ and $\Sym(\mu)$ is the number of $\s$ which leave
$\mu$ fixed. (Or equivalently, we take the sum over all distinct
permutations of the numbers $m_i$.) For example:
\[
    S_2(m)=R_2(m,0^4)+R_2(0,m,0^3)+R_2(0^2,m,0^2)+R_2(0^3,m,0)+R_2(0^4,m)
\]
\[
    =\frac{2m+1}{2m+3}Q_2(2m+3,1^4)
    +\sum_{q=0}^{2k-1}\frac13Q_2(1^{2k-q-1},m,1^q)
\]
\[
    =\frac{2m+7}5-\frac6{(2m+5)(2m+3)}
\]
If $\ll$ is any partition of $m$ then Theorem \ref{thm:recursive
formula for blln} says that
\[
    b_{\ll,k}^{m+k}=\sum_\mu \frac{b_\ll^\mu
    S_k(\mu)}{(-2)^{k+1}(2k-1)!!}
\]
where the sum is over all partitions $\mu$ of $m$ with at most
$2k+1$ parts. This gives a recursive formula for $b_\ll^m$. The
coefficients $b_\ll^\mu$ are then given by the sum of products
formula (Lemma~\ref{lem:sum of products formula for bllmu}).

\end{document}

%% file: gtoutput.tex
\def\ifplaintex{\expandafter\ifx\csname documentclass\endcsname\relax}

\ifplaintex 
\else
\fi

\expandafter\ifx\csname beginpicture\endcsname\relax
\expandafter\ifx\csname documentclass\endcsname\relax
\input pictex \else
\input prepictex \input pictex \input postpictex \fi\fi

\def\gt{{\mathsurround=0pt\it $\cal G\mskip-2mu$eometry \&\ 
$\cal T\!\!$opology}}        

\def\gtp{{\mathsurround=0pt\it $\cal G\mskip-2mu$eometry \&\ 
$\cal T\!\!$opology $\cal P\!$ublications}}  

\def\lognumber#1{\def\thelognumber{#1}}
\def\volumenumber#1{\def\thevolumenumber{#1}}
\def\papernumber#1{\def\thepapernumber{#1}}
\def\volumeyear#1{\def\thevolumeyear{#1}}

\def\pagenumbers#1#2{\def\startpage{#1}\def\finishpage{#2}}
\def\published#1{\def\publishdate{#1}}
\def\proposed#1{\def\theproposer{#1}}
\def\seconded#1{\def\theseconders{#1}}
\def\received#1{\def\receiveddate{#1}}

\def\accepted#1{\def\accepteddate{#1}}

\long\def\asciiabstract#1{\long\def\theasciiabstract{#1}}
\def\asciikeywords#1{\def\theasciikeywords{#1}}

\let\\\par\let\thelognumber\relax
\let\thevolumenumber\relax\let\thepapernumber\relax
\let\thevolumeyear\relax\let\thesamplenumber\relax\let\startpage\relax
\let\finishpage\relax\let\publishdate\relax\let\receiveddate\relax
\let\reviseddate\relax\let\accepteddate\relax\let\theasciititle\relax
\let\theasciiauthors\relax
\let\theasciiabstract\relax\let\theasciikeywords\relax
\let\theasciiemail\relax\let\theshortauthors\relax\let\theshorttitle\relax

\long\def\maketitlep{   

\count0=\startpage

\gt\hfill      
\beginpicture
\setcoordinatesystem units <0.33truein, 0.33truein> point at 2.2 0.9
\setplotsymbol ({$\cal G$})
\plotsymbolspacing=9truept
\circulararc 315 degrees from 0 1 center at 0 0
\setplotsymbol ({$\cal T$})
\circulararc 315 degrees from 1 -1 center at 1 0
\endpicture
\break
{\small\ifx\thesamplenumber\relax 
Volume \else Sample
\fi\thevolumenumber\ (\thevolumeyear)
\startpage--\finishpage\nl
Published: \publishdate}
\vglue 0.5truein plus 0.4fil minus 0.1truein

{\parskip=0pt\leftskip 0pt plus 1fil\def\\{\par\smallskip}{\ifplaintex\large
\else\Large\fi\bf\thetitle}\par\medskip}   

\vglue 0pt plus 0.1fil 

{\parskip=0pt\leftskip 0pt plus 1fil\def\\{\par}{\sc\theauthors}
\par\medskip}

\vglue 0pt plus 0.1fil 

{\small\parskip=0pt\let\newline\\
{\leftskip 0pt plus 1fil\def\\{\par}{\sl\theaddress}\par}
\expandafter\ifx\theemail\relax    
\relax\else\vglue 5pt plus 0.02fil minus 2pt\def\\{\stdspace{\rm 
and}\stdspace} 
\cl{Email:\stdspace\tt\theemail}\fi
\ifx\theurl\relax                  
\relax\else\vglue 5pt plus 0.02fil minus 2pt\def\\{\stdspace{\rm 
and}\stdspace}
\cl{URL:\stdspace\tt\theurl}\fi\par}

\vglue 7pt plus 0.3fil minus 3pt

{\bf Abstract}
\vglue 5pt plus 0.1fil minus 2pt

\theabstract

\vglue 7pt plus 0.3fil minus 3pt

{\bf AMS Classification numbers}\quad Primary:\quad \theprimaryclass

Secondary:\quad \thesecondaryclass

\vglue 5pt plus 0.3fil minus 2pt

{\bf Keywords:}\quad \thekeywords

\vglue 10pt plus 0.5fil minus 5pt

{\small  Proposed: \theproposer\hfill Received: \receiveddate\nl
Seconded: \theseconders\hfill 
\ifx\reviseddate\relax                         
Accepted: \accepteddate                        
\else
Revised: \reviseddate                          
\fi}
\eject
}       

\let\maketitlepage\maketitlep
\let\maketitle\maketitlepage

\font\phead=cmsl9 scaled 950
\font=cmsl9 scaled 1050
\font\pnum=cmbx10 scaled 913
\font\lnum=cmbx10 
\font\pfoot=cmsl9 scaled 950
\font=cmsl9 scaled 1050
\ifplaintex
\headline{\vbox to 0pt{\vskip -4.5mm\line{\small\phead\ifnum
\count0=\startpage ISSN 1364-0380 (on line)
1465-3060 (printed) \hfill {\pnum\folio}\else\ifodd\count0\def\\{ }%
\ifx\theshorttitle\relax\thetitle\else\theshorttitle\fi\hfill{\pnum\folio}
\else\def\\{ and }{\pnum\folio}\hfill\ifx\theshortauthors\relax\theauthors
\else\theshortauthors\fi\fi\fi}\vss}}
\footline{\vbox to 0pt{\vglue 0mm\line{\small\pfoot\ifnum\count0=\startpage
\copyright\ \gtp\hfill\else
\gt, Volume \thevolumenumber\ (\thevolumeyear)\hfill\fi}\vss
}}
\else
\makeatletter
\def\@oddhead{{\small\lhead\ifnum\count0=\startpage ISSN 1364-0380 (on line)
1465-3060 (printed) \hfill {\lnum\number\count0}\else\ifodd\count0
\def\\{ }\ifx\theshorttitle\relax \thetitle \else\theshorttitle\fi\hfill
{\lnum\number\count0}\else\def\\{ and }{\lnum\number\count0}
\hfill\ifx\theshortauthors\relax 
\theauthors\else\theshortauthors\fi\fi\fi}}\def\@evenhead{\@oddhead}
\def\@oddfoot{\small\lfoot\ifnum\count0=\startpage\copyright\ \gtp\hfill\else
\gt, Volume \thevolumenumber\ (\thevolumeyear)\hfill\fi}
\def\@evenfoot{\@oddfoot}
\makeatother
\fi

\newwrite\gtoutfile
\long\gdef\makeheadfile{  
{\def\\{, }\def\s{ }
\immediate\openout\gtoutfile head.xxx
\immediate\write\gtoutfile{Proxy-for: \ifx\theasciiauthors\relax
\theauthors\else\theasciiauthors\fi\s<\ifx\theasciiemail\relax\theemail\else\theasciiemail\fi>}
\immediate\write\gtoutfile{\noexpand\\}
\immediate\write\gtoutfile{Authors: \ifx\theasciiauthors\relax
\theauthors\else\theasciiauthors\fi}
{\def\\{ }\immediate\write\gtoutfile{Title: \ifx\theasciititle\relax
\thetitle\else\theasciititle\fi}}
\immediate\write\gtoutfile{Subj-class: GT or SG or MG etc}
\immediate\write\gtoutfile{MSC-class: \theprimaryclass\ifx\thesecondaryclass\relax\else, \thesecondaryclass\fi}
\immediate\write\gtoutfile{Journal-ref: Geom. Topol. \thevolumenumber
(\thevolumeyear) \startpage-\finishpage}
\immediate\write\gtoutfile{Comments: Published by Geometry and Topology at}
\immediate\write\gtoutfile{\s\s http://www.maths.warwick.ac.uk/gt/GTVol\thevolumenumber/paper\thepapernumber.abs.html}
\immediate\write\gtoutfile{\noexpand\\}
\immediate\write\gtoutfile{}
\ifx\theasciiabstract\relax
\immediate\write\gtoutfile{\theabstract}\else
\immediate\write\gtoutfile{\theasciiabstract}\fi
\immediate\write\gtoutfile{}
\immediate\write\gtoutfile{\noexpand\\}
\immediate\write\gtoutfile{}
\immediate\closeout\gtoutfile}}  

\def\maketitlepage{\maketitlep\makeheadfile}
\let\maketitle\maketitlepage

%% file: 2004-26.bbl
\begin{thebibliography}
\bibitem{[Arbarello-Cornalba:96]}
{\bf E Arbarello}, {\bf M Cornalba}, {\it Combinatorial and algebro-geometric
  cohomology classes on the moduli spaces of curves}, J. Alg. Geom. {5}
  (1996) 705--749

\bibitem{[Culler-Vogtmann-86]}
{\bf Marc Culler}, {\bf Karen Vogtmann}, {\it Moduli of graphs and automorphisms of free groups}, Invent. Math. {84} (1986)  91--119

\bibitem{[CV02]}
{\bf James Conant}, {\bf Karen Vogtmann}, {\it On a theorem of {K}ontsevich},
  \agtref3{2003}{42}{1167}{1224}

\bibitem{[I:MMM_and_Witten]} {\bf Kiyoshi Igusa}, {\it Combinatorial
{M}iller--{M}orita--{M}umford classes and {W}itten cycles},
\agtref4{2004}{23}{473}{520}

\bibitem{[I:BookOne]} {\bf Kiyoshi Igusa}, {\it Higher
{F}ranz--{R}eidemeister {T}orsion}, AMS/IP Studies in Advance
Mathematics, vol.~31, International Press (2002)

\bibitem{[I:GraphCoh]} {\bf Kiyoshi Igusa}, {\it Graph cohomology and
{K}ontsevich cycles}, to appear in Topology, \arxiv{math.AT/0303157}

\bibitem{[I:ComplexTorsion]}
{\bf Kiyoshi Igusa}, {\it Higher complex torsion and the framing principle}, 
to appear in Memoirs of AMS  \arxiv{math.KT/0303047}

\bibitem{[Kontsevich:Airy]} {\bf Maxim Kontsevich}, {\it Intersection
theory on the moduli space of curves and the matrix {A}iry function},
Comm. Math. Phys. {147} (1992) 1--23

\bibitem{[Miller86:MMM]}
{\bf Edward~Y. Miller}, {\it The homology of the mapping class group}, J.
  Differential Geom. {24} (1986) 1--14

\bibitem{[Mondello]}
{\bf Gabriele Mondello}, {\it Combinatorial classes on the moduli space of curves
  are tautological}, to appear in IMRN, \arxiv{math.AT/0303207}

\bibitem{[Morita84]} {\bf Shigeyuki Morita}, {\it Characteristic
classes of surface bundles}, Bull. Amer.  Math. Soc.  {11}
(1984) 386--388

\bibitem{A=B}
{\bf Marko Petkov\v{s}ek}, {\bf Herbert~S Wilf}, {\bf Doron Zeilberger}, {\it {$A=B$}}, published by A\,K Peters Ltd. Wellesley, MA (1996)

\bibitem{StanleyEC1} {\bf Richard~P Stanley}, {\it Enumerative
combinatorics. {V}ol. 1}, Cambridge Studies in Advanced Mathematics,
vol.~49, Cambridge University Press, Cambridge (1997) with a foreword
by Gian-Carlo Rota (corrected reprint of the 1986 original)

\bibitem{[Strebel]} {\bf Kurt Strebel}, {\it Quadratic
{D}ifferentials}, Springer-Verlag, Berlin (1984)

\end{thebibliography}
